\documentclass[12pt,reqno]{article}
\usepackage{amsfonts, amsthm, amsmath, amssymb, amscd}

\theoremstyle{plain}
\newtheorem{theorem}{Theorem}[section]
\newtheorem{lemma}[theorem]{Lemma}
\newtheorem{proposition}[theorem]{Proposition}
\newtheorem{corollary}[theorem]{Corollary}

\theoremstyle{definition}

\theoremstyle{remark}
\newtheorem{remark}[theorem]{Remark}

\numberwithin{theorem}{section}

\numberwithin{equation}{section}

\newcommand \lap {\lambda^{\prime}}
\newcommand \pip {\pi^{\prime}}

\newcommand \la {\lambda}
\newcommand \gz {{\cal G}}

\newcommand \R {{\cal R}}
\newcommand \A {{\cal A}}
\newcommand \W {{\cal W}}
\newcommand \T {{\cal T}}
\newcommand \q {{\bf q}}
\newcommand \ba {{\bf a}}

\newcommand \p {{\bf p}}

\newcommand \U {{\cal U}}

\newcommand \Y {{\cal Y}}

\newcommand \F {{\cal F}}

\newcommand \wab {{\cal W}_{{\cal A},\,B}}

\newcommand \omabz {{\Omega}_{{\cal A},B}^{{\mathbb Z}}}
\newcommand \omab { {\Omega}_{{\cal A},B}}

\newcommand \omzq {{\Omega}_{\q}^{{\mathbb Z}}}

\newcommand {\modk} {{\cal M}_{\kappa}}

\newcommand \HH {{\cal H}}

\begin{document}
\title{ Existence and Uniqueness of the Measure of Maximal Entropy for the
Teichm{\"u}ller Flow on the Moduli Space of Abelian Differentials.}

\author{Alexander I. Bufetov\footnote{Department of Mathematics,
Rice University, and the Steklov Institute of Mathematics, Russian
Academy of Sciences.} \ and Boris M. Gurevich\footnote{Department of
Mechanics and Mathematics, Moscow State University, and the
Institute for Information Transmission Problems, Russian Academy of
Sciences.}}


\date{}

\maketitle

\section{Introduction}
\label{introd} The Teichm\" uller geodesic flow $\{g_t\}$, first
studied by H. Masur \cite{masur} and W. Veech \cite{veech}, acts on
the moduli space of Riemann surfaces endowed with a holomorphic
differential.  More precisely, let $S$ be a closed surface of genus
$g\ge2$. One introduces on $S$ a complex structure $\sigma$ and a
holomorphic differential $\omega$. The pair $(\sigma,\omega)$ is
considered to be equivalent to another pair of the same nature
$(\sigma_1,\omega_1)$ if there is a diffeomorphism of $S$ sending
$(\sigma,\omega)$ to $(\sigma_1,\omega_1)$. The moduli space
$\mathcal M(g)$ consists of the equivalence classes, and the flow
$\{g_t\}$ on $\mathcal M(g)$ is induced by the action on the pairs
$(\sigma,\omega)$ defined by the formula
$g_t(\sigma,\omega)=(\sigma',\omega')$, where
$\omega'=e^t{\Re}(\omega)+ie^{-t}{\Im}(\omega)$, while the complex
structure $\sigma'$ is determined by the requirement that $\omega'$
be holomorphic. If $(\sigma,\omega)$ and $(\sigma',\omega')$ are
equivalent, then the differentials $\omega$ and $\omega'$ have the
same orders of zeros and the same area. Therefore, these orders and
area are well-defined on $\mathcal M(g)$. Moreover, they are
preserved by the Teichm\" uller flow $\{g_t\}$. Take an arbitrary
non-ordered collection $\kappa=(k_1,\dots, k_r)$ with $k_i\in\mathbb
N$, $k_1+\dots+k_r=2g-2$, and denote by $\mathcal M_\kappa$ the
subspace of $\mathcal M(g)$ corresponding to the differentials of
area 1 (i.e., $(i/2)\int\omega\wedge\bar\omega=1$) with orders of
zeros $k_i,\,i=1,\dots,r$; $\mathcal M_\kappa$ is said to be a {\it
stratum} in $\mathcal M(g)$. Each stratum is a $\{g_t\}$-invariant
set, and there is a natural $\{g_t\}$-invariant measure on $\mathcal
M_\kappa$; this measure is finite \cite{masur}, \cite{veech}.

The space $\mathcal M_\kappa$ also admits a natural topological
structure, in which it is in general non-connected. The number of
connected components is no more than 3 and depends on $\kappa$ (see
\cite{KZ} for details), each of them is $\{g_t\}$-invariant.

We fix an arbitrary closed component $\mathcal H$ and denote by
$\mu_\kappa$ the normalized restriction to $\mathcal H$ of the
above-mentioned $\{g_t\}$-invariant measure.

Veech \cite{veechteich} showed that $\{g_t\}$ with respect to the
measure $\mu_{\kappa}$ is a Kolmogorov flow with entropy given by
the formula
\begin{equation}
\label{entropy} h_{\mu_{\kappa}}(\{g_t\})=2g-1+r.
\end{equation}
Our aim is to establish the following
\begin{theorem}
\label{maxentropy} The measure $\mu_{\kappa}$ is the unique measure
of maximal entropy for the flow $\{g_t\}$ on $\mathcal H$.
\end{theorem}

The proof of this theorem is based on the representation of the flow
$\{g_t\}$ as a suspension flow over a countable alphabet topological
Markov shift. The reasoning proceeds in two steps. We begin with
sufficient conditions for an invariant measure of the above
suspension flow to be a measure with maximal entropy. These
conditions are contained in Theorem \ref{main} stated in Subsection
\ref{over_markov}. After stating the theorem we outline its proof.
This proof is close in spirit to thermodynamic formalism for
countable alphabet topological Markov shifts \cite{GS}, \cite{S},
\cite{sarig}. In particular, we use a uniqueness theorem by Buzzi
and Sarig \cite{BS} for an equilibrium measure. An application of
thermodynamic formalism  to another smooth dynamical system with
non-compact phase space, the geodesic flow on the modular surface,
can be found in \cite{GK}.

Subsections \ref{susp_entr}--\ref{proof2} are devoted to the proof
of Theorem \ref{main} in detail. In the rest of the paper (Sections
\ref{flow}--\ref{abel}) we deduce Theorem \ref{maxentropy} from
Theorem \ref{main}, and there (at the beginning of Section
\ref{flow}), as before, we start  from a sketch of the subsequent
reasoning.

The following observation lies at the centre of our argument in this
part of the proof. The Teichm{\"u}ller flow admits infinitely smooth stable and
unstable foliations, with respect to
which it is ``measurably Anosov'' in the sense of Veech
\cite{veechteich} and Forni \cite{forni}. The Masur-Veech measure
$\mu_{\kappa}$ induces globally defined sigma-finite measures on
unstable leaves and these measures are uniformly expanded by the
flow. In other words, the Masur-Veech measure has the Margulis
\cite{margulis} uniform expansion property on unstable leaves.
Informally,  Proposition \ref{unifexp} expresses the Margulis
property in terms of the symbolic representation of the
Teichm{\"u}ller flow.

Let us remark that to establish Theorem \ref{maxentropy} we need
only the special case of Theorem \ref{main} dealing with the
countable alphabet topological Bernoulli shift. But the proof for
this case would be only a little easier than in the general one. The
main results of this paper are stated without proof in \cite{BG}.

\section {Suspension flows}
\label{suspens}

Let $G$ be an Abelian group (in what follows only $G=\mathbb Z$ or
$G=\mathbb R$ will appear), and let $\{T_g,\,g\in G\}$ be an action
of $G$ by measurable transformations of a metrizable topological
space $X$ endowed with its Borel $\sigma$-algebra $\mathcal B$. Two
actions, $\{X,T_g\}=\{T_g,\,g\in G\}$ on $(X,\mathcal B)$ and
$\{X',T'_g\}=\{T'_g,\,g\in G\}$ on $(X',\mathcal B')$, are called
{\it isomorphic} if there is a one-to-one epimorphic bimeasurable
map $\Phi:X\to X'$ such that $T'_g\circ\Phi=\Phi\circ T_g$ for all
$g\in G$. If $\Phi$ is not necessarily epimorphic, we say that
$\{X,T_g\}$ is {\it embedded} into $\{X',T'_g\}$.

Consider also an action $\{X,T_g\}$ together with a
$\{T_g\}$-invariant Borel probability measure $\mu$ on $X$. Denote
such an object by $\{X,T_g;\mu\}$. We say that $\{X,T_g;\mu\}$ and
$\{X',T'_g;\mu'\}$ are isomorphic if there are sets $X_1\in\mathcal
B$, $X'_1\in\mathcal B'$ invariant with respect to all $T_g$ and all
$T'_g$ respectively such that $\mu(X_1)=\mu'(X'_1)=1$ and the
restrictions $\{X_1,T_g|_{X_1}\}$ and $\{X_1',T'_g|_{X'_1}\}$ are
isomorphic in the above sense.

If $G=\mathbb Z$, the corresponding action will be denoted by
$\{X,T_n\}$. If $G=\mathbb R$, we write $\{X,T_t\}$ or $\{T_t\}$ (or
just $T_t$, $S_t$, etc. when it cannot cause confusion). In the
former case $T_n=T^n$ where $T$ is a bimeasurable one-to-one
transformation of $(X,\mathcal B)$ called an {\it automorphism}. In
the latter case the action is called a {\it flow}. We keep the same
terminology for actions considered together with their invariant
measures. In this paper we mostly deal with flows that can be
defined as follows. Let $T$ be an automorphism of $(X,\mathcal B)$
and $f:X\to[c,\infty)$, $c>0$, be a measurable function. Consider
the direct product $X\times\mathbb R_+$ and its subspace
$X_f=\{(x,u):x\in X,\,0\le u<f(x)\}$. For $t\ge0$ and every point
$\tilde x=(x,u)\in X_f$, we set $S_t\tilde x=(x,u+t)$ if $u+t<f(x)$,
and $S_t\tilde x=(T^nx,u+t-\sum_{i=0}^{n-1}f(T^ix))$, where $n$ is
such that $\sum_{i=0}^{n-1}f(T^ix)\le u+t<\sum_{i=0}^nf(T^ix)$. For
$t<0$ we set $S_t=(S_{-t})^{-1}$ and thus obtain a flow $\{X_
f,S_t\}$. For this flow, we shall also use the notation $(T,f)$ and
call it the {\it suspension flow} constructed by $T$ and the {\it
roof function} $f$.

Denote by $\mathcal M_{T,f}$ the set of all $T$-invariant Borel
probability measures $\mu$ on $X$ with $\mu(f)<\infty$ (here and in
the sequel, $\mu(f):=\int fd\mu$). Every $(T,f)$-invariant Borel
probability measure $\mu_f$ on $X_f$ is induced by a measure
$\mu\in\mathcal M_{T,f}$. Namely,
$$\mu_f=(\mu(f))^{-1}(\mu\times\lambda)|_{X_f},
$$
where $\lambda$ is the Lebesgue measure on $\mathbb R_+$. We will
refer to $\mu_f$ as the $f$-{\it lifting} of $\mu$ and for brevity
write $(T,f;\mu_f)$ instead of $((T,f);\mu_f)$. The entropy
$h(T,f;\mu_f)$ of the flow $(T,f)$ with respect of the measure
$\mu_f$ is given by Abramov's formula
\begin{equation}
\label{abramov} h(T,f;\mu_f)=h(T;\mu)/\mu(f),
\end{equation}
where $h(T;\mu)$ is the entropy of the automorphism $T$ with respect
to the measure $\mu$. We define the {\it topological entropy} of
$(T,f)$ by
\begin{equation}
\label{topent} h_{\text{top}}(T,f)=\sup_{\mu\in\mathcal
M_{T,f}}h(T,f;\mu_f).
\end{equation}
This terminology is justified be the following well-known fact: if
$X$ is a compact space, $T$ is a homeomorphism of $X$, and $f$ is
continuous, then the right-hand side of (\ref{topent}) is indeed the
topological entropy of the suspension flow $(T,f)$.

We refer to every $\mu\in\mathcal M_{T,f}$ at which the supremum in
(\ref{topent}) is achieved as to a {\it measure of maximal entropy}
for $(T,f)$.

\subsection{Suspension flows over Markov shifts}
\label{over_markov}

In the specific case we will deal with, $(X,T)$ is a countable
alphabet topological Markov shift, i.e., $X$ is the set of infinite
two-sided paths of a directed graph $\Gamma=(V,E)$ with vertex set
$V$ and edge set $E\subseteq V\times V$, and $T$ is the shift
transformation: $(Tx)_i=x_{i+1}$ for each $x=(x_i,\,i\in\mathbb
Z)\in X$. In other words, $X$ consists of all sequences $x\in
V^{\mathbb Z}$ such that $B_{x_i,x_{i+1}}=1$, where $B=B(\Gamma)$ be
the incidence matrix of the graph $\Gamma$. The vertices $v\in V$
will also be called letters.

In the sequel we assume that $\Gamma$ is connected. If $\Gamma$ is
the complete graph, i.e., $E=V\times V$, we have the {\it
topological Bernoulli shift with alphabet $V$}.

We introduce the discrete topology on $V$, the product topology on
$V^\mathbb Z$, and the induced topology on $X\subset V^\mathbb Z$.
The map $T$ is clearly a homeomorphism of $X$. We shall refer to
every finite path of $\Gamma$, i.e., a sequence
$w=(v_1,\dots,v_k)\in V^k$ such that $(v_i,v_{i+1})$,
$i=1,\dots,k-1$, as a {\it word} and sometimes say that this word is
from $v_1$ to $v_k$. Denote the set of all words (including the
empty word) by $W(\Gamma)$ .

Let $w=(v_1,\dots,v_k)$, $w'=(v'_1,\dots,v'_l)$ be two words. The
concatenation $ww':= (v_1,\dots,v_k,v'_1,\dots,v'_l)$ is also a word
if $(v_k,v'_1)\in E$. We say that $w$ contains $w'$ (or $w'$ is a
subword of $w$) if $v'_1=v_i,\dots,v'_l=v_{i+l-1}$ for some $i$,
$1\le i\le k-l+1$. In the special case that $i=1$, we call $w'$ a
{\it prefix} of $w$. Let $w=(v_1,\dots,v_n)$, $n\ge 2$, be a word
and $w'=(v_1,\dots,v_l)$, $l\le n$, be a prefix of $w$. We call $w'$
a {\it simple prefix} of $w$ if there is no $k$, $2\le k\le l$, such
that $(v_1,\dots,v_{n-k+1})=(v_k,\dots,v_n)$. If $w$ is a simple
prefix of itself, then $w$ is called a {\it simple word}. If a
simple word is a prefix of another word, it is clearly a simple
prefix.
\begin{remark}
\label{words} Every word $w=(v_1,\dots,v_n)$ is certainly the
concatenation of the single-letter words $v_i$, so we will also
write $w=v_1\dots v_n$.
\end{remark}

To every word $w$ we assign the cylinder $C_w=\{x \in
X:(x_0,\dots,x_{|w|-1})=w\}$, where $|w|$ is the length of $w$,
i.e., the number of symbols in $w$.

For a function $f:X\to\mathbb R$, we set
\begin{equation*}
\label{var} \text{var}_n(f)=\sup\{|f(x)-f(y)|:x_i=y_i\ \text{when}\
|i|\le n\},\ \ n\in\mathbb N.
\end{equation*}
We say that $f$ has {\it summable variations} if
$\sum_{n=1}^\infty\text{var}_n(f)<\infty$, and that $f$ {\it depends
only on the future} if $x_i=y_i$ for all $i\ge 0$ implies that
$f(x)=f(y)$.

For a suspension flow $\{S_t\}=(T,f)$ and for a set $C\subset X$, we
put
\begin{equation}
\label{tau}\tilde\tau(x,C)=\inf\{t>0:S_t(x,0)\in C\times\{0\}\},\ \
x\in X,
\end{equation}
so that $\tilde\tau(x,C)$ is the first hitting time of
$C\times\{0\}$ for a point $x\times\{0\}\in X_f$.
\begin{theorem}
\label{main} For a countable alphabet topological Markov shift
$(X,T)$ corresponding to a connected graph $\Gamma$, let
$f:X\to[c,\infty)$, $c>0$, be a function with summable variations
depending only on the future, and let $\{S_t\}=(T,f)$ be the
suspension flow constructed by $T$ and $f$. Assume that
$\bar\mu\in\mathcal M_{T,f}$ is a measure positive on all cylinders
in $X$, and that for each $l>0$, there exists a simple word $w\in
W(\Gamma)$ with $|w|>l$ such that for every word $\hat w$ that does
not contain $w$ and for $\bar\mu$-almost all $x\in C_{w\hat ww}$, we
have
\begin{equation}
\label{marg} |\bar\mu(C_{w\hat
ww})/\bar\mu(C_w)-e^{-s\tilde\tau(x,C_w)}|\le
e^{-\alpha|w|-s\tilde\tau(x,C_w)},
\end{equation}
where $\tilde\tau(x,C_w)$ is defined in (\ref{tau}) and $\alpha,\,s$
are positive constants (depending only on $\bar\mu$). Then

(i) $s=h_{\text{top}}(T,f)$,

(ii) if $s=h(T,f;\bar\mu_f)$, then $\bar\mu_f$ is the unique measure
of maximal entropy for the flow $\{S_t\}=(T,f)$.
\end{theorem}
\begin{remark}
\label{cond_prob} The ratio in the left hand side of \eqref{marg} is
clearly the conditional measure of $C_{w\hat ww}$, given $C_w$.
\end{remark}
\begin{remark}
\label{future} The assumption that $f$ depends only on the future is
made just for convenience: Theorem \ref{main} remains true without
this assumption, but in the sequel we use it only in the above
particular form.
\end{remark}

Let us outline the proof of Theorem \ref{main}. At the first stage
we consider (in Subsection \ref{susp_entr}) the particular case
where $(X,T)$ is a Bernoulli shift, while $f(x)$, $x\in X$, depends
only on $x_0$ (we then say that $f$ depends on the zeroth
coordinate). In this case the topological entropy $h_{top}(T,f)$ can
be expressed explicitly in terms of $f$. At the next stage we come
back to the general case and prove that the supremum in the
definition of $h_{top}(T,f)$ can be taken over ergodic measures that
are positive on all cylinders in $X$ (see Subsection
\ref{positivemeas}). This enables us to state that
$h_{top}(T,f)=h_{top}(T_C,f_C)$ for the suspension flow $(T_C,f_C)$
where $T_C$ is the transformation induced by $T$ on a cylinder
$C\subset X$ and $f_C$ is determined naturally by $f$, $T$, and $C$.
If we chose $C=\{x:(x_0,\dots,x_k)=w\}$, where $w\in W(\Gamma)$,
then $(C,T_C)$ is isomorphic to the countable alphabet Bernoulli
shift. Hence the flow $(T_C,f_C)$ is isomorphic to a suspension flow
$(\sigma,\varphi)$ built over this Bernoulli shift. Here we use a
finite approximation and find a function $\varphi^w$ that depends on
the zeroth coordinate and is uniformly close to $\varphi$ (when $w$
is long enough). The topological entropies of the suspension flows
$(\sigma,\varphi)$ and $(\sigma,\varphi^w)$ are also close to each
other. We apply the results obtained at the first stage to the
latter flow and rewrite inequality (\ref{marg}) for it (see
Subsection \ref{proof2}). By directing the length of $w$ to infinity
we complete the proof of the equality $s=h_{top}(T,f)$, which
implies that the measure under consideration has maximal entropy. We
reduce the uniqueness of such a measure to that of the corresponding
equilibrium measure and here use the uniqueness theorem by Buzzi and
Sarig \cite{BS}.

\subsection{Entropy of suspension flows over Bernoulli shifts}
\label{susp_entr}

The proof of Theorem \ref{main} is based essentially on some
properties of suspension flows constructed by a topological Markov
shift (in particular, by a Bernoulli shift) and functions of one or
finitely many coordinates. Some of these properties, studied first
by Savchenko \cite{S}, are described in this section. We include
proofs for the reader's convenience. Our approach is close to that
of \cite{S}.

We begin with two simple lemmas. Let $\mathcal N=\mathbb N$ or
$\{1,\dots,n\}$, $n\ge 2$, and let $\mathbf c=(c_i,i\in\mathcal N)$
be a sequence of real numbers such that $\inf_{i\in\mathcal
N}c_i>0$. Denote by $\mathcal P=\mathcal P_{\mathcal N,\mathbf c}$
the family of sequences $\mathbf p=(p_i,i\in\mathcal N)$ such that
\begin{equation}
\label{p-seq} p_i\ge 0\ (i\in \mathcal N),\ \ \sum_{i\in\mathcal
N}p_i=1,\ \ \sum_{i\in\mathcal N}p_ic_i<\infty.
\end{equation}
(Certainly, $\mathcal P_{\mathcal N,\mathbf c}$ does not depend on
$\mathbf c$ when $|\mathcal N|<\infty$.) Let
\begin{equation}
\label{H-func} H(\mathbf p)=H_{\mathcal N,c}(\mathbf
p):=-\left(\sum_{i\in\mathcal N}p_i\log p_i\right)\left(\sum_{i\in
N}p_ic_i\right)^{-1},\ \ \mathbf p\in\mathcal P
\end{equation}
(we as usual let $0\log 0=0$).
\begin{lemma}
\label{posit} If $\mathbf p\in\mathcal P$ is such that $p_k=0$ for
some $k\in\mathcal N$, then there exists $\mathbf
p'=(p'_i,i\in\mathcal N)\in\mathcal P$ with $p'_i>0$ for all $i$
such that $H(\mathbf p')\ge H(\mathbf p)$, where the inequality is
strict when $H(\mathbf p)<\infty$.
\end{lemma}
\begin{proof}
We divide $\mathcal N$ into two non-empty subsets, $\mathcal
N^0=\{i\in\mathcal N:p_i=0\}$ and $\mathcal N^1=\mathcal
N\setminus\mathcal N^0$. Fix an arbitrary $l\in\mathcal N^1$ and for
$t\in[0,p_l)$ let $\mathbf p^t=(p_i^t,i\in\mathcal N)$, where
$p_k^t=t$, $p_l^t=p_l-t$, and $p_i^t=p_i$ for $i\ne k,l$. (By
assumption, $k\in\mathcal N^0$.) Clearly, $\mathbf p^t\in\mathcal P$
and $H(\mathbf p^t)=\infty$ when $H(\mathbf p)=\infty$.  A simple
calculation shows that if $H(\mathbf p)<\infty$, then the right-hand
derivative $\frac{d^+}{dt}H(\mathbf p^t)$ at $t=0$ is $+\infty$.
Hence $H(\mathbf p^t)>H(\mathbf p)$ when $t>0$ is small enough.

If $\mathcal N^0=\{k\}$, the proof is completed. If $\mathcal
N^0\setminus\{k\}\ne\emptyset$, we first consider the case
$H(\mathbf p)<\infty$. Fix an arbitrary $t\in(0,p_l)$ for which
$H(\mathbf p^t)>H(\mathbf p)$. It is easy to find positive numbers
$q_i$, $i\in\mathcal N^0\setminus\{k\}$, such that
\begin{equation}
\label{q-cond} \sum_{i\in\mathcal N^0\setminus\{k\}}q_i=1,\ \
\sum_{i\in\mathcal N^0\setminus\{k\}}q_i(c_i-\log q_i)<\infty.
\end{equation}
For $s\in[0,t)$ we put $\mathbf p^{t,s}=(p_i^{t,s},i\in\mathcal N)$,
where
\begin{equation}
\label{ts-pert} p_i^{t,s}=sq_i,\,i\in\mathcal N^0\setminus\{k\};\ \
p_i^{t,s}=p_i,\,i\in\mathcal N^1\setminus\{l\};\ \ p_k^{t,s}=t-s;\ \
p_l^{t,s}=p_l-t.
\end{equation}
From (\ref{q-cond}) it follows that $\mathbf p^{t,s}\in\mathcal P$
and $\lim_{s\to 0}H(\mathbf p^{t,s})=H(\mathbf p^t)$. Therefore,
$H(\mathbf p^{t,s})>H(\mathbf p)$ as $s>0$ is small enough, and
since $p_i^{t,s}>0$, we can take $\mathbf p^{t,s}$ with one of these
$s$ for $\mathbf p'$.

It remains to note that if $H(\mathbf p)=\infty$, then $H(\mathbf
p^{t,s})=\infty$ for all $s\in[0,t)$ (see (\ref{q-cond}),
(\ref{ts-pert})).
\end{proof}
\begin{lemma}
\label{suprem} Let $\mathcal N=\{1,\dots,n\}$, $n\ge 2$, and let
$\mathbf c$, $\mathcal P=\mathcal P_{\mathcal N}$, $H=H_{\mathcal
N,\mathbf c}$ be as above. Then $\sup_{\mathbf p\in\mathcal
P}H(\mathbf p)$ is the unique solution to the equation
$F_n(\beta)=1$, where $F_n(\beta)=\sum_{i=1}^ne^{-\beta c_i}$.
\end{lemma}
\begin{proof}
Since $H$ is a continuous function on the compact set $\mathcal
P\subset\mathbb R^n$, its supremum is attained at a point $\mathbf
p^0=(p_i^0,\,i=1,\dots,n)\in\mathcal P$. By Lemma \ref{posit}
$p_i^0>0$ for all $i$. Let $\mathcal P^+=\{\mathbf p\in\mathcal
P:p_i>0,\,i=1,\dots,n\}$. For $\mathbf p\in\mathcal P^+$ we put
$p_1=1-\sum_{i=2}^np_i$ and consider the equations $\partial
H(\mathbf p)/\partial p_i=0$, $i\in\mathcal N\setminus\{1\}$. From
this system we derive that if $\mathbf p^0$ is a point of extremum
of $H(\mathbf p)$, then $p_i^0=e^{-\beta^0c_i}/F_n(\beta^0)$, $1\le
i\le n$, where $\beta^0=\text{const}>0$. Hence the statement we are
proving is true when $c_i=c_1$ for $i=2,\dots,n$. Otherwise we take
any $i$ for which $c_i\ne c_1$ and from the equation $\frac{\partial
H(\mathbf p)}{\partial p_i}|_{\mathbf p=\mathbf p^0}=0$, where
$p_1=1-\sum_{i=2}^n$, obtain $\beta^0=H(\mathbf p^0)$. On the other
hand, by substituting $\mathbf p^0$ for $\mathbf p$ in $H(\mathbf
p)$ we see that $H(\mathbf
p^0)=\beta^0-\frac{F_n(\beta^0)}{F'_n(\beta^0)}\log F_n(\beta^0)$.
Therefore, $\log F_n(\beta^0)=0$, i.e., $\beta^0$ is a root of the
equation $F_n(\beta)=1$. This root is unique, since $F_n(\beta)$
decreases in $\beta$. Finally, $H(\mathbf p^0)=\max_{\mathbf
p\in\mathcal P}H(\mathbf p)$, because, as was mentioned above, every
point of maximum belongs to $\mathcal P^+$, hence the equations
$\partial H(\mathbf p)/\partial p_i=0$, $i\in\mathcal
N\setminus\{1\}$, $p_1=1-\sum_{i=2}^np_i$ must hold at this point.
But we already know that these equations have only one solution.
\end{proof}

Let us now consider a countable alphabet topological Bernoulli shift
$(X,T)$ with $X=V^\mathbb Z$, and the suspension flow
$\{S_t\}=(T,f)$ constructed by $T$ and a function $f$ such that
$f(x)=f_0(x_0)$, $x=(x_i,\,i\in\mathbb Z)$, where
$f_0:V\to[c,\infty)$, $c>0$. Let
$$F(\beta)=\sum_{v\in V} e^{-\beta f_0(v)},\ \ \beta\ge0.$$
\begin{lemma}
\label{b-flow} If there exists $\beta_0\ge0$ with $F(\beta_0)=1$,
then $h_{\text{top}}(T,f)=\beta_0$. Otherwise
$h_{\text{top}}(T,f)=\sup\{\beta\ge0:F(\beta)=\infty\}$.
\end{lemma}
\begin{proof}
Denote by $B_{T,f}$ the family of all Bernoulli measures in
$\mathcal M_{T,f}$. Each $\nu\in B_{T,f}$ is determined by the
one-dimensional distribution $\{p^\nu(v),\,v\in V\}$, where
$$p^\nu(v)=\nu(C_v)\ge 0,\ \ \sum_{v\in V}p^\nu(v)=1,\ \
\sum_{v\in V}p^\nu(v)f_0(v)<\infty.
$$
We note that
\begin{equation}
\label{supsup} \sup_{\mu\in\mathcal M_{T,f}}[h(T;\mu)/\mu(f)]=
\sup_{\nu\in B_{T,f}}[h(T;\mu)/\mu(f)].
\end{equation}

Indeed, every $\mu\in\mathcal M_{T,f}$ gives rise to the measure
$\mu_B\in B_{T,f}$ with $p^{\mu_B}(v)=\mu_B(C_v)=\mu(C_v)$. Clearly,
$\mu_B(f)=\mu(f)$, and basic properties of the measure--theoretic
entropy imply that $h(T;\mu)\le h(T;\mu_B)$.

Let us number in an arbitrary way the elements $v\in V$ and put
$B^{(n)}=\{\nu\in B_{T,f}:p^\nu(v_i)=0\ \text{for}\ i\ge n+1\}$,
$n\in\mathbb N$. For each $\mu\in B_{T,f}$, one can easily find a
sequence of measures $\nu_n\in B^{(n)}$ such that
$$\lim_{n\to\infty}[h(T;\nu_n)/\nu_n(f)]=h(T;\mu)/\mu(f).$$
Therefore, by (\ref{supsup}),
\begin{equation}
\label{supfin} \sup_{\mu\in\mathcal M_{T,f}}[h(T;\mu)/\mu(f)]=
\sup_{n\in\mathbb N}\sup_{\nu\in B^{(n)}}[h(T;\nu)/\nu(f)].
\end{equation}

We now notice that the relations
\begin{equation*}
\label{bernentropy} p_i:=p^\nu(v_i),\ 1\le i\le n;\ \ \mathbf
p=\mathbf p^\nu:=(p_1,\dots,p_n)
\end{equation*}
establish a one-to-one correspondence between $B^{(n)}$ and
$\mathcal P=\mathcal P_\mathcal N$ with $\mathcal N=\{1,\dots,n\}$,
and that $h(T;\nu)/\nu(f)=H_{\mathcal N,\mathbf c}(\mathbf
p)=H(\mathbf p)$, where $\mathbf c=(c_i,\,i\in\mathcal N)$,
$c_i=f_0(v_i)$, $i\in\mathcal N$ (see (\ref{p-seq}) and
(\ref{H-func})).

By Lemma \ref{suprem} the right-hand side of (\ref{supfin}) is
$\sup_n\beta_n$, where $\beta_n$ is determined by $F_n(\beta_n)=1$.
Let us note that $F_n$ is the $n$th partial sum of the series for
$F$ and that both $F_n$ and $F$ are strictly decreasing functions
(for $F$ it is true on the semi-axis where $F$ is finite). Hence
$\sup_n\beta_n=\lim_{n\to\infty}\beta_n$. We consider two possible
cases and first suppose that $F(\beta)=\infty$ for all $\beta\ge0$.
It is clear that in this case $\lim_{n\to\infty}\beta_n=\infty$.
Otherwise there exists a unique $\beta_\infty>0$ such that either
$F(\beta_\infty)=1$, or $F(\beta)<1$ for $\beta\ge\beta_\infty$ and
$F(\beta)=\infty$ for $\beta<\beta_\infty$. Since
$F_n(\beta)<F_{n+1}(\beta)<F(\beta)$ for all $n\ge1$ and
$\beta\ge0$, in both cases we have
$\lim_{n\to\infty}\beta_n\le\beta_\infty$. If
$\lim_{n\to\infty}\beta_n=:\beta'_{\infty}<\beta_\infty$, then
$F(\beta'_n)>1$ (in the latter case $F(\beta'_\infty)=\infty$).
Therefore $F_n(\beta'_\infty)>1$ for $n$ large enough. But
$\beta'_\infty>\beta_n$, hence $F_n(\beta'_\infty)<F_n(\beta_n)<1$
for all $n$. From this we conclude that
$\lim_{n\to\infty}\beta_n=\beta_\infty$. We thus come to both
statements of the lemma.
\end{proof}

\subsection{Induced automorphisms and Markov-Bernoulli reduction}
\label{induced_aut}

For the next lemma we have to remind the following definition. Let
$T$ be an automorphism of the space $(X,\mathcal B)$, and let
$C\in\mathcal B$. Denote
\begin{equation}
\label{def-c'} X_C=\{x\in X:\sum_{n<0}\mathbf
1_C(T^nx)=\sum_{n>0}\mathbf 1_C(T^nx)=\infty\}, \ \ C'=C\cap X_C.
\end{equation}
Thus $C'$ consists of all points in $C$ that return to $C$
infinitely often in forward and backward time. Let also
\begin{equation}
\label{return} \tau(T,C;x)=\min\{n>0:T^nx\in C\},\ \
T_{C'}x=T^{\tau(T,C;x)}x,\ \ x\in C'.
\end{equation}
It is clear that the sets $X_C,\,C'$ are measurable and invariant
with respect to $T$ and $T_{C'}$ respectively, and that $T_{C'}$ is
an automorphism of the set $C'$ provided with the induced Borel
$\sigma$-algebra; $T_{C'}$ is said to be the {\it induced
automorphism} on $C'$.
\begin{lemma}
\label{topol} Let $(T,f)$ be the suspension flow constructed by an
automorphism $T$ of $(X,\mathcal B)$ and a $\mathcal B$-measurable
function $f:X\to [c,\infty)$, $c>0$, and let $C\in\mathcal B$. Then
the suspension flow $(T|_{X_C},f|_{X_C})$ constructed by the
restrictions of $T$ and $f$ to $X_C$ is isomorphic to the suspension
flow $(T_{C'},f_{C'})$, where
\begin{equation}
\label{func} f_{C'}(x)=\sum_{i=0}^{\tau(T,C;x)-1}f(T^ix),\ \ x\in
C'.
\end{equation}
Furthermore, if $\mu\in\mathcal M_{T,f}$ is ergodic and such that
$\mu(C)>0$, then $\mu(f)=\int_{C'}f_{C'}d\mu$ and the suspension
flow $(T,f;\mu_f)$ is isomorphic to the suspension flow
$(T_{C'},f_{C'};(\mu_{C'})_{f_{C'}})$, where $\mu_f$ is the
$f$-lifting of $\mu$, $\mu_{C'}$ is the normalized restriction of
$\mu$ to $C'$, and $(\mu_{C'})_{f_{C'}}$ is the $f_{C'}$-lifting of
$\mu_{C'}$.
\end{lemma}

We omit the proof of this lemma, since it follows immediately from
standard facts of ergodic theory (see, for instance, \cite{CFS}).

The following construction is reminiscent of the Doeblin first
return method in the theory of Markov chains and has appeared repeatedly
in the literature in different forms (presumably for the first time
--- in \cite{G0}, see also \cite{G} and \cite{GS}).

Let $w=(v_1,\dots,v_l)\in W(\Gamma)$ and $C=C_w$. Then $X_C$ defined
by (\ref{def-c'}) can be described as follows: $x\in X$ belongs to
$X_C$ if and only if there is an increasing sequence of integers
$i_k=i_k(x)$, $-\infty<k<\infty$, such that $i_k\le 0$ for $k\le 0$,
$i_k>0$ for $k>0$, and $(x_{i_k},\dots,x_{i_k+l-1})=w$ for every
$k$, while no other segment of $x$ agrees with $w$. Furthermore,
$C'$ consists of those $x$ for which $i_0(x)=0$. It is clear that
\begin{equation*}
\label{return1} i_1(x)=\tau(T,C';x),\ \ i_k(x)\ge i_1(x)+k-1,\ \
x\in C'.
\end{equation*}
Denote by $A_w$ the set of all words $w'=(v'_1,\dots,v'_{l'})\in
W(\Gamma)$ with $l'>l$ such that
$(v'_1,\dots,v'_{l})=(v'_{l'-l+1},\dots,v'_{l'})=w$ and no other
subword of $w'$ (i.e., a word of the form
$(v'_m,v'_{m+1},\dots,v'_n)$, $1\le m\le n\le l'$) agrees with $w$.
It is easy to see that if $x\in X_C$, then for each $k\in\mathbb Z$,
the word $(x_{i_k},x_{i_k+1},\dots,x_{i_{k+1}+l-1})$ belongs to
$A_w$. We thus obtain a mapping $\Psi_w:X_C\to(A_w)^\mathbb Z$
measurable with respect to the appropriate Borel $\sigma$-algebras;
its restriction to $C'$ obviously induces a one-to-one
correspondence between $C'$ and $(A_w)^\mathbb Z$. Moreover, if
$x\in C'$, then $\Psi_wT_{C'}x=\sigma_w\Psi_wx$, where $\sigma_w$ is
the shift transformation on $Y_w:=(A_w)^\mathbb Z$, i.e., $(\sigma_w
y)_i=y_{i+1}$, $y=(y_i,\,i\in\mathbb Z)\in Y_w$. Therefore, $T_{C'}$
is isomorphic to the countable alphabet Bernoulli shift
$(Y_w,\sigma_w)$ with alphabet $A_w$. Here and in the sequel we
consider each $a\in A_w$ as either a word in the alphabet $V$ or a
letter in the new alphabet $A_w$. What of these two possibilities
takes place will always be clear from the context.

This construction reduces in essence the study of the topological
Markov shift $(X,T)$ to that of a topological Bernoulli shift
determined by $w$, and so we shall refer to it as the {\it
Markov--Bernoulli (M--B) reduction} applied to $(X,T)$ and $w$.

\subsection{Positive measures}
\label{positivemeas}

Our next aim is to show that the topological entropy of a suspension
flow over a Markov shift can be computed using only ergodic measures
that are positive on all cylinders.
\begin{lemma}
\label{sup} Let $(X,T)$ and $f$ be as in Theorem \ref{main} and let
$(T,f)$ be the suspension flow constructed by $T$ and $f$. Then
$$h_{\text{top}}(T,f)=\sup_{\mu\in\mathcal E^+_{T,f}}h(T,f;\mu_f),
$$
where $\mathcal E^+_{T,f}$ consists of all ergodic measures in
$\mathcal M_{T,f}$ that are positive on all cylinders in $X$.
\end{lemma}
\begin{proof}
Denote by $\mathcal E_{T,f}$ the set of ergodic measures in
$\mathcal M_{T,f}$. If $\mu\in\mathcal M_{T,f}\setminus\mathcal
E_{T,f}$, i.e., if $\mu$ is non-ergodic with respect to $T$, then
$\mu_f$, the $f$-lifting of $\mu$, is non-ergodic with respect to
the suspension flow $(T,f)$. The flow $(T,f;\mu_f)$ can be
decomposed into ergodic components (see \cite{R1}). This means the
following. There exists a measurable partition $\zeta$ of the space
$(X_f,\mu_f)$ such that $\mu_f$-almost every element $C_\zeta$ of
$\zeta$ is $(T,f)$-invariant and the conditional measure
$(\mu_f)^{C_\zeta}$ induced by $\mu$ on $C_\zeta$ is invariant and
ergodic with respect to the restriction of $(T,f)$ to $C_\zeta$. We
may consider $(\mu_f)^{C_\zeta}$ as a measure on the whole space
$X_f$; it is $(T,f)$-invariant and ergodic. By a general formula
(\cite{R2}, Section 9)
$$
h(T,f;\mu_f)=\int_{X_f|\zeta}h(T,f;(\mu_f)^{C_\zeta})
\mu_{f,\zeta}(dC_\zeta),
$$
where $\mu_{f,\zeta}$ is the projection of $\mu_f$ on the quotient
space $X_f|\zeta$. Hence, for every $\varepsilon>0$, there is an
element $C_\zeta$ with
$h(T,f;(\mu_f)^{C_\zeta})>h(T,f;\mu_f)-\varepsilon$. On the other
hand, $(\mu_f)^{C_\zeta}$, being a $(T,f)$-invariant probability
measure on $X_f$, is the $f$-lifting of a $T$-invariant probability
measure $\mu^{C_\zeta}$ on $X$, i.e.,
$(\mu_f)^{C_\zeta}=(\mu^{C_\zeta})_f$. It is clear that
$\mu^{C_\zeta}(f)<\infty$ and $(T,\mu^{C_\zeta})$ is ergodic. Since
$\varepsilon>0$ was arbitrary small, we conclude that

\begin{equation*}
\label{topent2} h_{\text{top}}(T,f)=\sup_{\mu\in\mathcal
E_{T,f}}h(T,f;\mu_f).
\end{equation*}

Let $\mathcal E_{T,f}^0:=\mathcal E_{T,f}\setminus\mathcal
E_{T,f}^+$ and assume that, contrary to the lemma we have to prove,
for some $\delta\in(0,\infty)$,
\begin{equation}
\label{bigger} \sup_{\mu\in\mathcal E_{T,f}^0}h(T;\mu)/\mu(f)>
\sup_{\mu\in\mathcal E_{T,f}^+}h(T;\mu)/\mu(f)+\delta,
\end{equation}
which in particular means that
$$
\sup_{\mu\in\mathcal
E_{T,f}^+}h(T;\mu)/\mu(f)<\infty.
$$
By virtue of (\ref{bigger})
there is $\mu^0\in\mathcal E_{T,f}^0$ such that
\begin{equation}
\label{ent1} h(T;\mu^0)/\mu^0(f)\ge\sup_{\mu\in\mathcal
E_{T,f}^+}h(T;\mu) /\mu(f)+\delta/2.
\end{equation}

To show that this is impossible we first consider the case
$h(T;\mu^0)<\infty$ and let $h^0=h(T,\mu^0)/\mu^0(f)$. Since $f$ has
summable variations, one can find $n_\delta\in\mathbb N$ such that,
for every $n\ge n_\delta$, there is a function $f_n:X\to\mathbb R_+$
with the following three properties: $f_n(x)=f_n(y)$ whenever
$x_i=y_i$ for $|i|\le n$, $\inf_{x\in X}f_n(x)\ge c$, and
$$
\sup_{x\in X}|f(x)-f_n(x)|\le\delta c^2/8h^0.
$$
One can easily check that then
\begin{equation}
\label{difent} \left|\frac{h(T;\mu)}{\mu(f_n)}-\frac{h(T;\mu)}
{\mu(f)}\right|<\delta
\end{equation}
for every $\mu\in\mathcal E_{T,f}^+\cup\{\mu^0\}$. Hence (see
(\ref{ent1}))
\begin{equation}
\label{ent2} h(T;\mu^0)/\mu^0(f_n)\ge\sup_{\mu\in\mathcal
E_{T,f_n}^+} h(T;\mu)/\mu(f_n)+\delta/4.
\end{equation}
Since $|f-f_n|<\text{const}$, the functions $f$ and $f_n$ are
integrable or not integrable with respect to a finite measure
simultaneously. Hence $\mathcal E_{T,f_n}=\mathcal E_{T,f}$ and
$\mathcal E_{T,f_n}^+=\mathcal E_{T,f}^+$.

If $h(T;\mu^0)=\infty$, then (\ref{ent2}) clearly holds as well.

Using the assumption $\mu^0\in\mathcal E_{T,f}^0$, we find a word
$w^0\in W(\Gamma)$ with $\mu^0(C_{w^0})=0$. Fix an arbitrary
$n^1\ge\max\{n_\delta,\,|w^0|\}$ and a word $w^1\in W(\Gamma)$ with
$|w^1|=n^1$, $\mu^0(w^1)>0$. Then we set $f^1:=f_{n^1}$,
$C:=C_{w^1}$, and apply the M--B reduction to $(X,T)$ and $w^1$. By
Lemma \ref{topol} the suspension flow $(T|_{X_C},f^1|_{X_C})$ is
isomorphic to the suspension flow
$(\sigma,\varphi):=(\sigma_{w^1},\varphi_{f^1,w^1})$, where
\begin{equation}
\label{phi} \varphi(y):=(f^1)_{C'}(\Psi_{w^1}^{-1}y), \ \ y\in
Y_{w^1}.
\end{equation}
Notice that the function $\varphi$ is constant on every
one-dimensional cylinder $\{y\in Y:y_0=a\}$, $a\in A_{w^1}$; the
reason is that each $a\in A_{w^1}$ when considered as a word from
$W(\Gamma)$ is not shorter than $w^1$.

Let us carry over the measure $\mu_{C'}^0$, the normalized
restriction of $\mu^0$ to $C'$ (where $C'$ is defined in
(\ref{def-c'})), to $Y$ via the mapping $\Psi_{w^1}$ to obtain a
Borel probability measure $\nu^0$ on $Y$. From the above-described
properties of $\Psi_{w^1}$ it follows that the automorphisms
$(T_{C'};(\mu^0)_{C'})$ and $(\sigma;\nu^0)$ are isomorphic and
hence, by Lemma \ref{topol}, the suspension flow
$(T,f^1;(\mu^0)_{f^1})$ is isomorphic to the suspension flow
$(\sigma,\varphi;(\nu^0)_\varphi)$, where $(\mu^0)_{f^1}$ and
$(\nu^0)_\varphi$ are the $f^1$-lifting of $\mu^0$ and the
$\varphi$-lifting of $\nu^0$, respectively. Therefore, by
(\ref{abramov}),
\begin{equation}
\label{ent3} h(T;\mu^0)/\mu^0(f^1)=h(\sigma;\nu^0)/\nu^0(\varphi).
\end{equation}

If we change $\nu^0$ for a $\sigma$-invariant Bernoulli measure
$\nu^1$ with the same one-dimensional distribution (i.e., with
$\nu^1(C_a)=\nu^0(C_a)$ for all $a\in A_{w^1}$, where $C_a=\{y\in
Y:y_0=a\})$, then the numerator on the right-hand side of
(\ref{ent3}) can only increase, while the denominator will not
change (since $\varphi$ is constant on every cylinder $C_a$, $a\in
A_{w^1}$).

From the definition of $\nu^0$ and $\nu^1$ it follows that
$\nu^0(C_{a^0})=\nu^1(C_{a^0})=0$ for some $a^0\in A_{w^1}$. Indeed,
let $w^1=(v_1^1,\dots,v_{l_1}^1)$. Since the graph $\Gamma$ is
connected, there exists a word $(v_1,\dots,v_r)\in W(\Gamma)$ with
$(v_1,\dots,v_{l_1})=w^1$, $(v_{r-l_0+1},\dots,v_r)=w^0$, where
$l_0=|w^0|$. Choose an arbitrary shortest word of this type and
denote it by $w'$. Similarly, let $w''$ be one of the shortest words
in which there are an initial subword and a terminal subword that
coincide with $w^0$ and $w^1$, respectively. From the assumption
that $\mu^0(w^0)=0$, $\mu^0(w^1)>0$, $|w^0|\le|w^1|$ it follows that
$w''=w^0\hat w$ where $\hat w$ can be of one of the following three
forms: (a) $\hat w=w^1$; (b) $\hat w=\hat w^1w^1$, $\hat w^1\in
W(\Gamma)$; (c) $\hat w=(v_k^1,\dots,v_{l_1}^1)$, $1<k\le l_1$, is a
terminal subword of $w^1$. Consider the word $w'\hat w$. One easily
checks that $w'\hat w\in A_{w^1}$. Moreover, $\mu^0(w'\hat w)=0$,
because $w'\hat w$ contains $w^0$ as a subword. Hence
$\mu^0_{C'_{w^1}}(C'_{w^1}\cap C_{w'\hat w})=0$. We can put
$a^0:=w'\hat w$. Since $\Psi_{w^1}(C'\cap C_{w'\hat w})=C_{a^0}$, we
have $\nu^1(C_{a^0})=\nu^0(C_{a^0})=0$.

We now want to perturb $\nu^1$ within the class of Bernoulli
measures on $Y$ in such a way as to obtain a measure for which the
right-hand side of (\ref{ent3}) is bigger than for $\nu^0$ and which
is positive on all cylinders.

Since $\varphi(y),\ y=(y_i,\,i\in\mathbb Z)\in Y$, depends solely on
$y_0$, we have $\varphi(y)=\varphi_0(y_0)$, where $\varphi_0$ is a
function on $A_{w^1}$.

Using Lemma \ref{posit}, we find a $\sigma$-invariant Bernoulli
measure $\nu^2$ on $Y$ such that if $h(\sigma;\nu^1)<\infty$, then
\begin{equation}
\label{ent4}
\frac{h(\sigma;\nu^2)}{\nu^2(\varphi)}>\frac{h(\sigma;\nu^1)}
{\nu^1(\varphi)} \ge\frac{h(\sigma;\nu^0)}{\nu^0(\varphi)},
\end{equation}
and if $h(\sigma;\nu^1)=\infty$, then $h(\sigma;\nu^2)=\infty$ as
well.

Apply the mapping $\Psi_{w^1}^{-1}$ to transfer the measure $\nu^2$
to $C'$ and denote the resulting measure by $\mu'$. The suspension
flow $(\sigma,\varphi;(\nu^2)_\varphi)$ is then isomorphic to the
suspension flow $(T_{C'},(f^1)_{C'};(\mu')_{(f^1)_{C'}})$. Let
$C'(n)=\{x\in C':\tau_{T,C}(x)\}$, $n=1,2,\dots$, and
$\mu'_n=\mu'|_{C'(n)}$ be the restriction of $\mu'$ to $C'(n)$
considered as a measure on $X$. Then the measure
$$\mu'':=\sum_{n=1}^\infty\sum_{k=0}^{n-1}T^k\circ\mu'_n$$
is concentrated on $X_C$ and $T$-invariant. By normalizing $\mu''$
we obtain a probability measure $\mu'''$. By Lemma \ref{topol} the
flows $(T|_{X_C};f^1|_{X_C})$ and $(T_{C'},(f^1)_{C'})$ are
isomorphic. Then the flow $(T,f_1;(\mu''')_{f^1})$ is isomorphic to
the flow $(T_{C'},(f^1)_{C'};(\mu')_{{f^1}_C'})$ and hence (see
above) to the flow $(\sigma,\varphi;(\nu^2)_\varphi)$. Therefore,
\begin{equation}
\label{ent5} h(T;\mu''')/\mu'''(f^1)=h(\sigma;\nu^2)/\nu^2(\varphi)>
h(T;\mu^0)/\mu^0(f^1)
\end{equation}
(see (\ref{ent3}), (\ref{ent4})). It is clear that
$\mu'''\in\mathcal E_{T,f}$. Moreover, $\mu'''\in\mathcal
E_{T,f}^+$. Otherwise we could apply to $\mu'''$ the procedure that
lead us to the measure $\nu^0$, starting from $\mu^0$. The resulting
measure would coincide with $\nu^{t,s}$, and there would be a letter
$a\in A_{w^1}$ with $\nu^{t,s}(C_a)=0$. But we know that this is
impossible. Thus (\ref{ent5}) contradicts (\ref{ent2}) with
$f_n=f^1$ and hence contradicts (\ref{bigger}).
\end{proof}
\begin{corollary}
\label{topent1} Let $\Gamma$, $(X,T)$, $f$ be as in Theorem
\ref{main}, $(Y_w,\sigma_w)$ be the topological Bernoulli shift
obtained from $(X,T)$ and $w$ by the M--B reduction, where $w\in
W(\Gamma)$, and let $\varphi_{f,w}$ be the function defined in
(\ref{phi}). Then the suspension flows $(T,f)$ and
$(\sigma_w,\varphi_{f,w})$ have the same topological entropy.
\end{corollary}
\begin{proof}
As before, we let $C=C_w$ an use the notation in
(\ref{def-c'})--(\ref{func}). From the definition of $\sigma_w$ and
$\varphi_{f,w}$ it follows immediately that the suspension flows
$(T_{C'},f_{C'})$ and $(\sigma_w,\varphi_{f,w})$ are isomorphic and
hence
$h_{\text{top}}(T_{C'},f_{C'})=h_{\text{top}}(\sigma_w,\varphi_{f,w})$.
Similarly, by virtue of Lemma \ref{topol},
$h_{\text{top}}(T_{C'},f_{C'})=h_{\text{top}}(T|_{X_C},f|_{X_C})$.
But $h_{\text{top}}(T|_{X_C},f|_{X_C})\le h_{\text{top}}(T,f)$,
because $X_C$ is a $T$-invariant subset of $X$. Hence
\begin{equation}
\label{topent3}
h_{\text{top}}(\sigma_w,\varphi_{f,w})=(T_{C'},f_{C'}) \le
h_{\text{top}}(T,f).
\end{equation}
On the other hand, by the same Lemma \ref{topol}
$h(T_{C'},f_{C'};\tilde\mu_{C'})=h(T,f;\tilde\mu)$ for every
$\mu\in\mathcal E_{T,f}^+$, where $\tilde\mu_{C'}$ is the
$f_{C'}$-lifting of the normalized restriction of $\mu$ to $C'$, and
$\tilde\mu$ is the $f$-lifting of $\mu$. The supremum in
$\mu\in\mathcal E_{T,f}^+$ of the left-hand side of the last
equality is clearly not bigger than $h(T_{C'},f_{C'})$, while by
Lemma \ref{sup} the supremum of the right-hand side is
$h_{\text{top}}(T,f)$. Hence $h_{\text{top}}(T_{C'},f_{C'})\ge
h_{\text{top}}(T,f)$, which together with (\ref{topent3}) yields
what we are proving.
\end{proof}

\subsection {\bf Proof of Theorem \ref{main}}
\label{proof2}

Fix notation as in Theorem \ref{main}. Let us also fix $n$ and for a
wile write $w$ and $C$ instead of $w_n$ and $C_{w_n}$, respectively.
Consider the sets $X_C$, $C'$, the induced transformation
$T_{C'}:C'\to C'$, and the function $f_{C'}$ (see
(\ref{def-c'})--(\ref{func})). Apply the M--B reduction to $(X,T)$
and $w$.

From (\ref{tau}), (\ref{func}) it is clear that
$\tilde\tau(T,C;x)=f_{C'}(x)$ for every $x\in C'$. Thus (\ref{marg})
can be rewritten in the form
\begin{equation}
\label{cond1} |\bar\mu(C_{w\hat ww})/\bar\mu(C)-e^{-sf_{C'}(x)}| \le
e^{-\alpha|w|-sf_{C'}(x)},
\end{equation}
which is true for $\bar\mu$-almost all $x\in C_{w\hat ww}$.

The simplicity of $w$ implies that each word $a\in A_w$ is of the
form $a=w\hat ww$, where $\hat w\in W(\Gamma)$ ($\hat w$ may be an
empty word if $ww\in W(\Gamma)$), and $\hat w$ does not contain $w$
as a subword.

By assumption, the measure $\bar\mu$ is positive on all cylinders
and $T$-invariant. Hence $\bar\mu(C)=\bar\mu(C')>0$, and we can
normalize $\bar\mu$ on $C'$ to obtain a $T_{C'}$-invariant
probability measure $\bar\mu'$. Its image $\nu':=(\Psi_w)_*\mu'_0$
is a probability measure $\nu'$ on $Y=Y_w$ invariant with respect
the shift transformation $\sigma$. From the definition of $\Psi_w$
it follows that, for $a:=w\hat ww\in A_w$,
$$
C_a:=\{y\in Y:y_0=a\}=\Psi_w C_{w\hat ww}
$$
and hence
\begin{equation}
\label{nu} \nu'(C_a)=\mu_0(C_{w\hat ww})/mu_0(C_w).
\end{equation}
Notice that $\nu'(C_a)>0$ for all $a\in A_w$.

Taking into account the relation between $f_{C'}$ and
$\varphi=\varphi_{f,w}$ (see (\ref{phi})) and using (\ref{cond1}),
(\ref{nu}), we obtain
\begin{equation}
\label{cond2} |\nu'(C_a)-e^{-s\varphi(y)}|\le
e^{-\alpha|w|-s\varphi(y)},\ \ a\in A_w,
\end{equation}
for $\nu'$-almost all $y\in C_a$.

Our next step is to approximate $\varphi$ by a function that is
constant on each cylinder $C_a$, $a\in A_w$.

Since $w$ is simple, we have $\tau(T,C';x)\ge|w|$ for all $x\in C'$
(see (\ref{return})). We say that $x^{(1)}=(x_i^{(1)},\,i\in\mathbb
Z)\in C'$ and $x^{(2)}=(x_i^{(2)},\,i\in\mathbb Z)\in C'$ are
equivalent ($x^{(1)}\thicksim x^{(2)}$) if
$\tau(T,C';x^{(1)})=\tau(T,C';x^{(2)})$ and $x_i^{(1)}=x_i^{(2)}$
for $0\le i\le\tau(T,C';x^{(1)})$. If $x^{(1)}\thicksim x^{(2)}$,
then (because $w$ is simple) $x_i^{(1)}=x_i^{(2)}$ for
$\tau(T,C';x^{(1)})\le i\le\tau(T,C';x^{(1)})+|w|-1$ as well, from
which we obtain (see (\ref{func}))
\begin{align}
\label{estvar1}
|f_{C'}(x^{(1)})-&f_{C'}(x^{(2)})|\le \notag \\
&\le\sum_{i=0}^{\tau(T,C';x)-1}|f(T^ix^{(1)})-f(T^ix^{(2)})|\le
\sum_{n=|w|}^\infty\text{var}_n(f).
\end{align}

Let
\begin{equation}
\label{approx} C^w(x):=\{x'\in C':x'\thicksim x\},\ \
f^w(x)=\inf_{x'\in C^w(x)}f_{C'}(x').
\end{equation}
It is easy to see that $C^w(x)$ is a cylinder and that these
cylinders constitute a partition of $C'$. Moreover, by virtue of
(\ref{approx}), (\ref{estvar1}) the function $f^w$ is constant on
each element of this partition and
\begin{equation*}
\label{estvar2} 0\le
f_{C'}(x)-f^w(x)\le\sum_{n=|w|}^\infty\text{var}_n(f),\ \ x\in C'.
\end{equation*}
Therefore,
\begin{equation}
\label{estvar3}
0\le\varphi(y)-\varphi^w(y)\le\sum_{n=|w|}^\infty\text{var}_n(f),\ \
y\in Y,
\end{equation}
where $\varphi^w(y):=f^w(\Psi_w^{-1}y)$ is constant on each cylinder
$C_a\subset Y$, $a\in A_w$ (here, by $\Psi_w^{-1}y$ we mean the
unique point $x\in C'$ such that $\Psi_wx=y$) and hence there is a
function $\varphi_0^w$ on $A_w$ such that
$\varphi^w(y)=\varphi_0^w(y_0)$.

With Lemma \ref{b-flow} in mind we will estimate the sum $\sum_{a\in
A_w}\exp[-s\varphi_0^w(a)]$. Let
\begin{equation}
\label{vard} \delta_w:=\sum_{n=|w|}^\infty\text{var}_n(f).
\end{equation}
Since $\nu'(C_a)>0$ for all $a\in A_w$, one can choose, for every
$a$, a point $y_a\in C_a$ such that (\ref{cond2}) holds for $y=y_a$.
Hence
$$\nu'(C_a)-\exp[-\alpha|w|-s\varphi(y_a)]\le
e^{-s\varphi(y_a)}\le\nu'(C_a)+ \exp[-\alpha|w|-s\varphi(y_a)],
$$
so that
$$\nu'(C_a)/(1+e^{-\alpha|w|})\le e^{-s\varphi(y_a)}\le
\nu'(C_a)/(1-e^{-\alpha|w|}),\ \ a\in A_w,$$
$$1/(1+e^{-\alpha|w|})\le\sum_{a\in A_w}e^{-s\varphi(y_a)}\le
1/(1-e^{-\alpha|w|}).
$$
From (\ref{estvar3}), (\ref{vard}) we obtain
\begin{align}
\label{cond3} \frac{1}{1+e^{-\alpha|w|}}\le&\sum_{a\in
A_w}e^{-s\varphi(y_a)}
\le\sum_{a\in A_w}e^{-s\varphi_0^w(a)}\notag\\
&=\sum_{a\in A_w}e^{-s\varphi^w(y_a)}\le
\frac{1}{1-e^{-\alpha|w|}}+\sum_{a\in
A_w}\left[e^{-s\varphi^w(y_a)}-
e^{-s\varphi(y_a)}\right]\notag\\
=\frac{1}{1-e^{-\alpha|w|}}+&\sum_{a\in A_w}e^{-s\varphi(y_a)}
\left[e^{s(\varphi(y_a)-\varphi^w(y_a))}-1\right]\le
\frac{e^{s\delta_w}}{1-e^{-\alpha|w|}}.
\end{align}

By assumption, we now can take $w=w_n$, where $|w_n|\to\infty$ as
$n\to\infty$. From (\ref{cond3}) it follows that
\begin{equation}
\label{cond5} \lim_{n\to\infty}\sum_{a\in
A_{w_n}}\exp[-s\varphi_0^{w_n}(a)]=1.
\end{equation}
Let
$$F_n(u):=\sum_{a\in A_{w_n}}\exp[-u\varphi_0^{w_n}(a)],\ \ n=1,2,\dots
$$
If, for a fixed $n$, there is a $u\in\mathbb R$ such that $F_n(u)=1$
(such $u$ can be only one), then we denote this $u$ by $u_n$.
Otherwise we put $u_n:=\sup\{u:F_n(u)=\infty\}$. Notice that $u_n\ge
0$ (because $F_n(0)=\infty$) and $u_n<\infty$ (because of
(\ref{cond3})). From the definition of $\varphi^{w_n}$ it follows
that $\inf_{y\in Y}\varphi^{w_n}(y)\to\infty$ as $n\to\infty$
(remind that $Y=Y_{w_n}$). Therefore, for every $\gamma>0$, we have
$d F_n(u)/du \to-\infty$ as $n\to\infty$ uniformly in $u$ on the set
$D_\gamma:=\{u:\gamma<F_n(u)<\infty\}$ (we mean the right-side
derivative if $u$ is the left endpoint of the interval $D_\gamma$).
Using this fact, it is easy to deduce from (\ref{cond5}) that
$u_n\to s$ as $n\to\infty$ (it would be sufficient to know that $d
F_n(u)/du<\text{const}<0$ on $D_\gamma$).

Let us now consider two isomorphic suspension flows,
$(\sigma,\varphi^{w_n})$ and $(T_{C'},f^{w_n})$. By Lemma
\ref{b-flow}
$$t_n=h_{\text{top}}(\sigma,\varphi^{w_n})=
h_{\text{top}}(T_{C'},f^{w_n}),\ \ n=1,2,\dots,
$$
where $C'=(C_{w_n})'$, and hence
$$s=\lim_{n\to\infty}h_{\text{top}}(T_{C'},f^{w_n}).
$$
From (\ref{estvar3}) and the evident bounds
$h_{\text{top}}(\sigma,\varphi^{w_n})\le 2s$ (as $n$ is large
enough), $\inf\,\varphi\ge c$, and $\inf\,\varphi^{w_n}\ge c$ we
obtain
$$|h_{\text{top}}(\sigma,\varphi^{w_n})-
h_{\text{top}}(\sigma,\varphi)|\le 2s\delta_n/c,$$ where
$\delta_n=\sum_{k=|w_n|}\text{var}_n(f)$ (cf. (\ref{difent})).
Therefore $s=h_{\text{top}}(\sigma,\varphi)$, and hence (see
Corollary \ref{topent1}) $s=h_{\text{top}}(T,f)$. So statement (i)
is proved.
\medskip

To prove statement (ii) suppose that $s=h_{\rm
top}(T,f;(\bar\mu)_f)$. Together with (i) this means that
$(\bar\mu)_f$ is a measure with maximal entropy for the suspension
flow $(T,f)$. We conclude that $(T,f)$ can have only one measure
with maximal entropy, using a result by Buzzi and Sarig \cite{BS},
as follows.

Let $(\bar\mu)_f$ be such a measure. Then by (\ref{topent})
$$\frac{h(T,\mu)}{\mu(f)}\le\frac{h(T,\bar\mu)}{\bar\mu(f)}=s,\ \
\mu\in\mathcal M_{T,f},
$$
where $s=h_{\text{top}}(T,f)$. Hence for every $\mu\in\mathcal
M_{T,f}$, we have $h(T,\mu)+\mu(g)\le 0$, where $g(x):=-sf(x)$,
$x\in X$, while $h(T,\bar\mu)+\bar\mu(g)=0$, so that the topological
pressure of $g$ is zero and $\bar\mu$ is a $g$-equilibrium measure.
Using the natural projection $\pi:V^\mathbb Z\to V^{\mathbb Z_+}$ we
let $X_+=\pi X$ and $f_+(x_+)=f(x)$ for $x_+\in X^+$ and any
$x\in\pi^{-1}x_+$ (by assumption, $f$ is constant on the set
$\pi^{-1}x_+$, so that $f(x)$ depends only on $x_+$). It is easily
checked that $\pi Tx=T_+\pi x$, $x\in X$, where $T_+$ is the shift
transformation on $X_+$, and, moreover, that $\pi$ induces a
one-to-one correspondence between $\mathcal M_{T,f}$ and $\mathcal
M_{T_+,f_+}$, the set of $T_+$-invariant probability measures
$\mu_+$ on $X_+$ with $\mu_+(f_+)<\infty$. Let $\mu_{0+}\in\mathcal
M_{T_+,f_+}$ correspond to $\mu_0$. Then $\mu_{0+}$ is a
$g_+$-equilibrium measure, where $g_+=-sf_+$. Notice that the
one-sided Markov shift $T_+$ is topologically transitive (because
the graph $\Gamma$ is connected), the topological pressure of $g_+$
is zero (because this is the case for $g$), and $\sup_{x_+\in
X_+}\,g_+(x_+)<0$ (because $\inf_{x\in X}\,f(x)>0$). Thus, by
Theorem 1.1 from \cite{BS}, there can be only one $g_+$-equilibrium
measure. So the proof of Theorem \ref{main} is completed.

\section{A covering flow}
\label{flow}

At this point we start deducing Theorem \ref{maxentropy} from
Theorem \ref{main}. The aim of the present section is to recall the
construction of a flow that can be viewed (see Section \ref{abel}
below) as "covering" for the Techm\" uller flow $\{g_t\}$. We shall
show (see corollary 5.2) that our problem can be reduced to a
similar problem for this covering flow (denoted by $\{P^t\}$).

We first remind some constructions due to Rauzy \cite{rauzy}, Veech
\cite{veech}, and Zorich \cite{zorich} (see also \cite{viana}).
Using these constructions we obtain (in Subsection \ref{zipper}) the
covering flow as a suspension flow over a measurable transformation
defined on a bounded Borel set in a Euclidian space of finite
dimension. We next consider (in Section \ref{symbolic}) a symbolic
representation of the flow $\{P^t\}$ and show that it is, up to an
isomorphism, a suspension flow over a countable alphabet topological
Markov shift (denote the alphabet by $\mathcal A$). Theorem
\ref{main} cannot yet be applied directly to this suspension flow,
since, if for no other reason than that its roof function is not
bounded away from zero and has non-summable variations. That is why
we change the base (Poincar{\'e} section) of our flow for its
appropriate cylinder subset (we in fact use a family of cylinders)
and thus perform a change-over to a new suspension representation
(going back to Veech). The new suspension flow appears to be built
over a countable alphabet Bernoulli shift whatever cylinder set we
take (this is a Markov--Bernoulli reduction as defined in Subsection
\ref{induced_aut}). Not all cylinders are useful for us here, but
only those corresponding to admissible ``positive'' words
$w\in\bigcup_{n=1}^\infty\mathcal A^n$ (there is a canonical way to
assign a matrix with non-negative integer entries to each $w$; if
all the entries are positive, we refer to $w$ as a positive word).
It we change the positivity for a stronger requirement that each
words $w$ involved in the construction has a simple positive prefix
(the definition is given in Subsection \ref{over_markov}), it will
be possible to prove (see Subsections \ref{roof} and \ref{transit})
that the roof function has summable variations (it is even H\"older
continuous), while the measure $\mu_\kappa$ (see Section
\ref{introd}) induces on the base of our suspension flow an
invariant measure satisfying the requirements imposed on the measure
$\bar\mu$ in Theorem \ref{main}. It remains to note that for each
$\{P^t\}$-invariant ergodic probability measure $\nu$ with positive
entropy, there is a sufficiently large collection of words $w$ with
$\nu(w)>0$ that have a simple positive prefix (see Lemma
\ref{positive}). We thus have everything necessary for applying
Theorem \ref{main}.

\subsection{Induction maps}
\label{R_V_Z}

Let $\pi$ be a permutation of $m$ symbols, which will always be
assumed irreducible in the sense that $\pi\{1,\dots
,k\}=\{1,\dots,k\}$ implies $k=m$. The Rauzy operations $a$ and $b$
are defined by the formulas
$$
a\pi(j)=\begin{cases}
\pi j,&\text{if $j\leq\pi^{-1}m$,}\\
\pi m,&\text{if $j=\pi^{-1}m+1$,}\\
\pi(j-1),&\text{if $\pi^{-1}m+1<j\le m$;}
\end{cases}
$$
$$
b\pi(j)=\begin{cases}
\pi j,&\text{if $\pi j\leq \pi m$,}\\
\pi j+1,&\text{if $\pi m<\pi j<m$,}\\
\pi m+1,&\text{ if $\pi j=m$.}
\end{cases}
$$

These operations preserve irreducibility. The {\it Rauzy class}
$\mathcal R(\pi)$ is defined as the set of all permutations that can
be obtained from $\pi$ by application of the transformation group
generated by $a$ and $b$. From now on we fix a Rauzy class $\R$ and
assume that it consists of irreducible permutations.

For $i,j=1,\dots,m$, denote by $E^{ij}$ the $m\times m$ matrix whose
$(i,j){\rm th}$ entry is $1$, while all others are zeros. Let $E$ be
the identity $m\times m$-matrix. Following Veech \cite{veech},
introduce the unimodular matrices
\begin{equation} \label{mat_a}
A(a,\pi)=\sum_{i=1}^{\pi^{-1}m}E^{ii}+E^{m,\pi^{-1}m+1}+
\sum_{i=\pi^{-1}m}^{m-1}E^{i,i+1},
\end{equation}
\begin{equation} \label{mat_b}
A(b,\pi)=E+E^{m,\pi^{-1}m}.
\end{equation}
For a vector $\la=(\la_1,\dots,\la_m)\in{\mathbb R}^m$, we write
$$
|\la|=\sum_{i=1}^m\la_i.
$$
Let
$$
\Delta_{m-1}=\{\la\in {\mathbb R}^m:|\la|=1,\ \la_i>0 \text{ for }
i=1,\dots,m\}.
$$

One can identify each pair $(\lambda,\pi)$,
$\lambda\in\Delta_{m-1}$, with the {\it interval exchange map} of
the interval $I:=[0,1)$ as follows. Divide $I$ into the
sub-intervals $I_k:=[\beta_{k-1},\beta_k)$, where $\beta_0=0$,
$\beta_k=\sum_{i=1}^k\lambda_i$, $1\le k\le m$, and then place the
intervals $I_k$ in $I$ in the following order (from left to write):
$I_{\pi^{-1}1},\dots, I_{\pi^{-1}m}$. We obtain a piecewise linear
transformation of $I$ that preserves the Lebesgue measure.

The space $\Delta(\R)$ of interval exchange maps corresponding to
$\R$ is defined by
$$
\Delta(\R)=\Delta_{m-1}\times\R.
$$
Denote
$$
\Delta_{\pi}^+=\{\la\in\Delta_{m-1}| \ \la_{\pi^{-1}m}>\la_m\},\ \
\Delta_{\pi}^-=\{\la\in\Delta_{m-1}| \ \la_m>\la_{\pi^{-1}m}\},
$$
$$
\Delta^+(\R)=\cup_{\pi\in{\cal R}}\{(\pi,\la)|\
\la\in\Delta_{\pi}^+\},$$
$$
\Delta^-(\R)=\cup_{\pi\in{\cal R}}\{(\pi,\la)|\
\la\in\Delta_{\pi}^-\},$$
$$
\Delta^\pm(\R)=\Delta^+(\R)\cup\Delta^-(\R).
$$
The {\it Rauzy-Veech induction map} ${\cal
T}:\Delta^\pm(\R)\to\Delta(\R)$ is defined as follows:
\begin{equation}
\label{te} {\cal T}(\la,\pi)=\begin{cases}
(\frac{A(a,\,\pi)^{-1}\la}{|A(a,\,\pi)^{-1}\la|},a\pi), &\text{if
$\la\in\Delta_\pi^+$,}\\
(\frac{A(b,\,\pi)^{-1}\la}{|A(b,\,\pi)^{-1}\la|},\,b\pi), &\text{if
$\la\in\Delta_\pi^-$} .
\end{cases}
\end{equation}

One can check that $\cal T(\la,\pi)$ is the interval exchange map
induced by $(\la,\pi)$ on the interval $[0,1-\gamma]$, where
$\gamma=\min(\la_m,\la_{\pi^{-1}m})$; this interval stretches to
unit length.

Denote
\begin{equation}
\label{delta+-} \Delta^\infty(\mathcal R)=\bigcap_{n\ge0}\mathcal
T^{-n} \Delta^\pm(\mathcal R).
\end{equation}
Every $\cal T$-invariant probability measure is concentrated on
$\Delta^\infty(\mathcal R)$. On the other hand, a natural Lebesgue
measure defined on $\Delta(\mathcal R)$, which is finite, but
non-invariant, is also concentrated on $\Delta^\infty(\mathcal R)$.
Veech \cite{veech} showed that $\mathcal T$ has an absolutely
continuous ergodic invariant measure on $\Delta(\R)$, which is,
however, infinite.

Following Zorich \cite{zorich}, for
$(\la,\pi)\in\Delta^\infty(\mathcal R)$, we set
\begin{equation}
\label{an} n(\la,\pi)=\begin{cases} \min \{k>0:{\cal
T}^k(\la,\pi)\in\Delta^-(\mathcal R)\},&\text{if
$\la\in\Delta_{\pi}^+$;}\\
\min \{k>0:{\cal T}^k(\la,\pi)\in\Delta^+(\mathcal R)\},&\text{if
$\la\in\Delta_{\pi}^-$.}
\end{cases}
\end{equation}

The {\it Rauzy-Veech-Zorich induction map} $\cal G$ is defined by
the formula
\begin{equation}
\label{ge} {\cal G}(\la,\pi)={\cal T}^{n(\la,\pi)}(\la,\pi),\ \
(\la,\pi)\in\Delta^\infty(\mathcal R).
\end{equation}
\begin{theorem}[\rm Zorich\,\cite{zorich}]
\label{zorichthm} The map ${\cal G}$ has an ergodic invariant
probability measure $\nu$ absolutely continuous with respect to the
Lebesgue measure on $\Delta^\infty(\R)$. The density $\rho$ of this
measure is of the form
$$\rho(\la,\pi)=P_{\pi}(\la)/Q_{\pi}(\la),\ \ \la=(\la_1,\dots,\la_m),$$
where $P_{\pi}$ and $Q_{\pi}$ are homogeneous polynomials with
non-negative coefficients.
\end{theorem}

\subsection{Zippered rectangles}
\label{zipper}

Here we briefly recall the construction of the Veech space of
zippered rectangles. We use the notation of \cite{bufetov}.

{\it Zippered rectangles} associated with the Rauzy class $\R$ are
triples $(\la,\pi,\delta)$, where
$\la=(\la_1,\dots,\la_m)\in{\mathbb R}^m$, $\la_i>0$, $\pi\in{\cal
R}$, $\delta=(\delta_1,\dots,\delta_m)\in{\mathbb R}^m$, and the
vector $\delta$ satisfies the following inequalities:
\begin{equation}
\label{deltaone} \delta_1+\dots+\delta_i\leq 0,\ \ i=1,\dots,m-1,
\end{equation}
\begin{equation}
\label{deltatwo} \delta_{\pi^{-1}\,1}+\dots+\delta_{\pi^{-1}\,i}\geq
0, \ \ i=1, \dots, m-1.
\end{equation}
The set of all vectors $\delta$ satisfying (\ref{deltaone}),
(\ref{deltatwo}) is a cone in ${\mathbb R}^m$; we denote it by
$K(\pi)$.

For a zippered rectangle $(\la,\pi,\delta)$, we set
\begin{equation}
\label{heigt} h_r:=-\sum_{i=1}^{r-1}\delta_i+\sum_{i=1}^{\pi
r-1}\delta_{\pi^{-1}\,i},
\end{equation}
\begin{equation}
\label{area} Area\,(\la,\pi,\delta):=\sum_{r=1}^m\la_rh_r.
\end{equation}
(Our convention is $\sum_{i=u}^v...=0$ when $u>v$.) By
\eqref{deltaone}, \eqref{deltatwo} $h_r\ge 0$ for all r, and if we
relate the set $\mathcal Z:=\bigcup_{r=1}^mI_r\times[0,h_r]$ (a
union of rectangles in $\mathbb R^2$) to every triple
$(\la,\pi,\delta)$, then $Area(\la,\pi,\delta)$ becomes merely the
Lebesgue measure (area) of $\mathcal Z$. By appropriate
identification of intervals in the boundaries of different
rectangles $I_r\times[0,h_r]$ one obtains a compact Riemannian
surface and a 1-form on it. This procedure (due to Veech
\cite{veech}) is thoroughly described in the literature (see, for
example, \cite{viana1}, \cite{zorich}). We do not use it as such and
so omit details.

Denote by ${\cal V}(\R)$ the space of all zippered rectangles
corresponding to the Rauzy class $\R$, i.e.,
$$
{\mathcal V}(\R)=\{(\la,\pi,\delta):\la\in{\mathbb
R}^m_+,\,\pi\in\R,\,\delta\in K(\pi)\}.
$$
Let also
$$
{\mathcal V}^+(\R)=\{(\la,\pi,\delta)\in{\mathcal
V}(\R):\la_{\pi^{-1}m}>\la_m\},
$$
$$
{\mathcal V}^-(\R)=\{(\la,\pi,\delta)\in{\mathcal
V}(\R):\la_{\pi^{-1}m}<\la_m\},
$$
$$
{\cal V}^\pm(\R)={\cal V}^+(\R)\cup{\cal V}^-(\R).
$$
Veech \cite{veech} introduced the flow $\{P^t\}$ acting on ${\cal
V}(\R)$ by the formula
$$
P^t(\la,\pi,\delta)=(e^{t}\la,\pi,e^{-t}\delta),
$$
and the map $\U:{\cal V}^\pm(\R)\to{\cal V}(\R)$, where
$$
{\U}(\la,\pi,\delta)=\begin{cases}(A(\pi,a)^{-1}\la,a\pi,A(\pi,
a)^{-1}\delta),
&\text{if $\la_{\pi^{-1}m}>\la_m$,}\\
(A(\pi,b)^{-1}\la,b\pi,A(\pi,b)^{-1}\delta), &\text{if
$\la_{\pi^{-1}m}<\la_m$.}
\end{cases}
$$
(The inclusion $\U{\cal V}^\pm(\R)\subset{\cal V}(\R)$ is not
evident and should be proved; this was done in \cite{veech}.) The
map $\U$ and the flow $\{P^t\}$ commute on ${\cal V}^\pm(\R)$ and
both preserve the measure determined on ${\cal V}(\R)$ by the volume
form $Vol=d\la_1\dots d\la_md\delta_1\dots d\delta_m$. They also
preserve the area of a zippered rectangle (see (\ref{area})) and
hence can be restricted to the set
$$
{\cal V}^{1,\pm}(\R):=\{(\la,\pi,\delta)\in{\cal V}^\pm(\R):
Area(\la,\pi,\delta)=1\}.
$$
The restriction of the volume form $Vol$ to ${\cal V}^{1,\pm}(\R)$
induces on this set a measure $\mu_\R$ which is invariant under
$\cal U$ and $\{P^t\}$.

For $(\la,\pi)\in\Delta(\R)$, denote
\begin{equation}
\label{taulambda}
\tau^0(\la,\pi)=:-\log(|\la|-\min(\la_m,\la_{\pi^{-1}m})).
\end{equation}
From (\ref{mat_a}), (\ref{mat_b}) it follows that if
$\la\in\Delta_\pi^+\cup\Delta_\pi^-$, then
\begin{equation}
\label{tau1} \tau^0(\la,\pi)=-\log|A^{-1}(c,\pi)\la|,
\end{equation}
where $c=a$ when $\la\in\Delta_\pi^+$, and $c=b$ when
$\la\in\Delta_\pi^-$.

Next denote
$$\Y_1(\R):=\{x=(\la,\pi,\delta)\in{\cal
V}(\R):|\la|=1,\ Area(\la,\pi,\delta)=1\},
$$
\begin{equation*}
\label{tau3}\tau(x):=\tau^0(\la,\pi)\text{ for
}x=(\la,\pi,\delta)\in\Y_1(\R),
\end{equation*}
\begin{equation}
\label{phase} {\cal V}_{1,\tau}(\R):=\bigcup_{x\in\Y_1(\R),\ 0\leq
t\leq \tau(x)}P^tx.
\end{equation}
Using the map $\U$, we are going to transfer the flow $P^t$ to the
set ${\cal V}_{1,\tau}(\R)$ (or, more precisely, to its proper
subset).

It is easy to check that $\U P^{\tau(x)}x\in\Y_1(\R)$ for every
$x\in\Y_1(\R)\cap{\cal V}^\pm(\R)$. Identifying the points
$P^{\tau(x)}x$ and $\U P^{\tau(x)}x$, we can continue the trajectory
of $x$ by some distance. But it can happen that $\U
P^{\tau(x)}x\notin{\mathcal V}^{\pm}(\R)$, so that we cannot proceed
this way.

To make $\{P^t\}$ well defined on an invariant set we have to
somewhat reduce the domain of $\U$. Let
$$
{\cal V}^{1,\pm}_{\ne}(\R):=\{(\la,\pi,\delta)\in {\cal
V}^{1,\pm}(\R):a_m(\delta)\ne 0\},
$$
$$
{\cal V}_\infty(\R):=\bigcap_{n\in\mathbb Z}\U^n{\cal
V}^{1,\pm}_{\ne}(\R).
$$
Clearly
$\U^n$ is well-defined on ${\cal
V}_\infty(\R)$ for all $n\in\mathbb Z$.

We now set
$$
\Y(\R):=\Y_1(\R)\cap\mathcal V_\infty(\R),\ \ {\tilde{\cal V}}(\R):=
{\cal V}_{1,\tau}(\R)\cap{\cal V}_\infty(\R).
$$
The above identification enables us to define on $\tilde{\cal
V}(\R)$ a natural flow, for which we retain the notation $\{P^t\}$.
(Although the bounded positive function $\tau$ is not separated from
zero, the flow $\{P^t\}$ is well defined.)

Let us also note that $\mu_{\R}({\cal V}_{1,\tau}(\R))>0$. By a
theorem proven independently by Veech \cite{veech} and Masur
\cite{masur} $\mu_{\R}({\cal V}_{1,\tau}(\R))<\infty$, and we in
what follows assume that the restriction of $\mu_{\R}$ to ${\cal
V}_{1,\tau}(\R)$ is normalized to have total mass $1$. Since clearly
$\mu_{\R}({\cal V}_{1,\tau}(\R)\setminus\tilde{\cal V}(\R))=0$, we
can consider $\mu_{\R}$ to be defined on $\tilde{\cal V}(\R)$. This
measure is ergodic with respect to the flow $\{P^t\}$.
\begin{remark}
\label{variab} The presentation here differs from the one in Veech
\cite{veech} by a linear change of variable: the Veech vector
parameters $h$ and $a$ of a zippered rectangle $(\la,\pi,\delta)$
are expressed in terms of $\pi$ and $\delta$ by \eqref{heigt} and
the equations
\begin{equation*}
\label{aa} a_r=a_r(\delta)=-\sum_{i=1}^{r}\delta_{i},\ \
r=1,\dots,m.
\end{equation*}
\end{remark}

Following Zorich \cite{zorich}, denote
$$
\Y^+(\R)=\{x=(\la,\pi,\delta)\in\Y(\R):
\la\in\Delta_\pi^+,\,a_m(\delta)<0\},
$$
$$
\Y^-(\R)=\{x=(\la,\pi,\delta)\in\Y(\R):
\la\in\Delta_\pi^-,\,a_m(\delta)>0\},
$$
$$
\Y^\pm(\R)=\Y^+(\R)\cup\Y^-(\R),
$$
and let ${\Y_\infty^{\pm}(\R)}$ be the set of all $x\in
\Y^{\pm}(\R)$ for which there exist infinitely many positive $t$ and
infinitely many negative $t$ such that $P^tx\in\Y^{\pm}(\R)$.

Take $x\in{\Y_\infty^{\pm}(\R)}$, $x=(\la,\pi,\delta)$, and let
$\F(x)$ be the first return of $x$ to the transversal $\Y^{\pm}(\R)$
under the flow $\{P^t\}$. The map $\F$ is an extension of the map
$\gz$ to the space of zippered rectangles:
\begin{equation}
\label{ef} {\rm if }\
\F(\la,\pi,\delta)=(\lap,\pip,\delta^{\prime}),\ {\rm then}\
(\lap,\pip)=\gz(\lap,\pip).
\end{equation}
Note that $\F$ is invertible on $\Y_\infty^{\pm}(\R)$ and if
$x\in\Y^+(\R)$ (respectively, $x\in\Y^-(\R)$), then
$\F(x)\in\Y^-(\R)$ (respectively, $\F(x)\in \Y^+(\R)$). Moreover,
$$
\Y_\infty^{\pm}(\R)=\bigcap_{n\in\mathbb Z}\F^n\Y^{\pm}(\R).
$$

If $x=(\la,\pi,\delta)\in\Y_\infty^{\pm}(\R)$ and $\mathcal
G(\la,\pi)=\T^n(\la,\pi)$, then, by (\ref{te}), (\ref{ge}), the
first return time of $x$ to $\Y^{\pm}(\R)$ under the flow $\{P^t\}$
is
\begin{equation}
\label{tau2} \tau(\la,\pi)+ \dots +\tau(\mathcal
T^{n-1}(\la,\pi))=-\log|A^{-1}(c,c^{n-1}\pi)\dots A^{-1}(c,\pi)\la|,
\end{equation}
where $c=a$ when $\la\in\Delta_\pi^+$, and $c=b$ when
$\la\in\Delta_\pi^-$.

We finish this section with considering a relationship between the
probability measure $\mu_\R$ mentioned above and the measure $\nu$
introduced in Theorem \ref{zorichthm}.

Denote by $\mu_\R^1$ the $\F$-invariant probability measure induced
by $\mu_\R$ on $\Y_\infty^\pm(\R)$. Remark that if
$(\la,\pi,\delta)\in\Y_\infty^\pm(\R)$, then
$(\la,\pi)\in\Delta^\infty$.
\begin{lemma}{\rm (cf. \cite{veech}, \cite{zorich})}
\label{meas_proj} Let
$\tilde\psi:\Y_\infty^\pm(\R)\to\Delta^\infty(\R)$ be the map
defined by $\tilde\psi(\la,\pi,\delta)=(\la,\pi)$. Then
$\tilde\psi_*\mu_\R^1=\nu$.
\end{lemma}
\begin{proof}
Note that there is a natural Lebesgue measure on each of the spaces
$\tilde{\mathcal V}(\R)$, $\Y_\infty^{\pm}(\R)$, and
$\Delta^\infty(\R)$. Since $\mu_\R$ is proportional to the Lebesgue
measure on $\tilde{\mathcal V}(\R)$, from the definition of
$\{P^t\}$ it follows that $\mu_\R^1$ is absolutely continuous with
respect to the Lebesgue measure on $\Y_\infty^{\pm}(\R)$. Let
$\nu^1:=\tilde\psi_*\mu_\R^1$. It is clear that $\nu^1$ is a
probability measure absolutely continuous with respect to
$mes_\Delta$, the Lebesgue measure on $\Delta^\infty(\R)$, while by
Theorem \ref{zorichthm} the probability measure $\nu$ is equivalent
to $mes_\Delta$ and ergodic with respect to $\mathcal G$. Therefore,
$\nu^1$=$\nu$.
\end{proof}

\section{Symbolic representation of the covering flow}
\label{symbolic}

In this section we construct suspension flows over symbolic Markov
shifts that will be of great importance in the study of the flow
$\{P^t\}$. Using \cite{veech}, \cite{zorich}, we begin with a brief
description of a symbolic model for the map $\gz$.

\subsection{Symbolic dynamics for the mapping $\mathcal G$}
\label{symbol1}

We shall only deal with the interval exchanges $(\la,\pi)$ from
$\Delta^\infty(\R)$ (see (\ref{delta+-})), so that all iterations of
the map $\gz$ are defined. Our notation follows \cite{bufetov}.

Consider the alphabet
$$
{\cal A}:=\{(c,n,\pi)|\ c=a\ {\rm or}\ b,\ n\in {\mathbb N},\
\pi\in\R\}.
$$

For $w_1=(c_1,n_1,\pi_1)\in\A$, $w_2=(c_2,n_2,\pi_2)\in\A$, we set
$$
B(w_1,w_2)=
\left\{\begin{aligned}1\ & {\rm if}\ (c_1)^{n_1}\pi_1=\pi_2,\,c_2\neq c_1,\\
0\ \ & {\rm otherwise}\\
\end{aligned}\right.
$$
and thus define a function $B:{\cal A}\times{\cal A}\to\{0,1\}$. In
other terms, we have a directed graph $\Gamma_{{\cal A},\,B}=(V,E)$,
where $V=\A$ and where $(w_1,w_2)\in E$ if and only if
$B(w_1,w_2)=1$. From the definition of the Rauzy class $\R$ in
Subsection \ref{R_V_Z} it follows that the graph $\Gamma_{{\cal
A},\,B}$ is connected.

Introduce the space of words
$$
{\cal W}_{{\cal A},\,B}=\{w=w_1\dots w_n|\ w_i\in{\cal A},\ B(w_i,
w_{i+1})=1\ {\rm for}\ i=1,\dots,n\}.
$$
It is convenient to include the empty word in $\wab$. We use the
notation of Sub-section \ref{over_markov}. In particular, for a word
$w\in{\cal W}_{{\cal A},\,B}$, we denote by $|w|$ its length, i.e.,
the number of symbols in it; given two words
$w^{(1)},w^{(2)}\in{\cal W}_{{\cal A},\,B}$, we denote by
$w^{(1)}w^{(2)}$ their concatenation. Note that the word
$w^{(1)}w^{(2)}$ need not belong to ${\cal W}_{{\cal A},\,B}$,
unless a compatibility condition is satisfied by the last letter of
$w^{(1)}$ and the first letter of $w^{(2)}$.

To each nonempty word $w\in{\cal W}_{{\cal A},\,B}$ we assign a {\it
renormalization matrix} $A(w)$ as follows. If $w$ is a single-letter
word, $w=(c,n,\pi)\in\A$, we set (see (\ref{mat_a}), (\ref{mat_b}))
\begin{equation}
\label{renorm0} A(w)=A(c,\pi)A(c,c\pi)\dots A(c,c^{n-1}\pi);
\end{equation}
for $w\in{\cal W}_{{\cal A},\,B}$, where $w=w_1\dots w_n$,
$w_i\in\A$, we set
\begin{equation}
\label{renorm1} A(w)=A(w_1)\dots A(w_n).
\end{equation}

Consider the sequence spaces
$$
\Omega_{\A, B}=\{\omega=(\omega_0,\omega_1,\dots)|\ \omega_n\in\A, \
B(\omega_n, \omega_{n+1})=1 \ {\rm for \ all}\ n\in {\mathbb Z_+}\},
$$
and
$$
\Omega_{\A, B}^{\mathbb
Z}=\{\omega=(\dots,\omega_{-1},\omega_0,\omega_1,\dots)|\
\omega_n\in\A, \ B(\omega_n, \omega_{n+1})=1 \ {\rm for \ all}\
n\in{\mathbb Z}\}.
$$
Denote by $\sigma$ the one step left shift on both these spaces.

Let us now describe the coding map. For every letter
$w=(c,n,\pi)\in\A$, we set
\begin{equation}
\label{Delta_w1} \Delta(w)=\begin{cases}
\Delta^+(\R)\cap\{(\la,\pi)\in\Delta^\infty(\R)|\
n(\la,\pi)\}\ &\text{if}\ c=a,\\
\Delta^-(\R)\cap\{(\la,\pi)\in\Delta^\infty(\R)|\ n(\la,\pi)\}\
&\text{if}\ c=b.
\end{cases}
\end{equation}
In other words, when $c=a$ (resp., $c=b$), $\Delta(w)$ consists of
all points $(\la,\pi)\in\Delta^\infty(\R)\cap\Delta^+(\R)$ (resp.,
$(\la,\pi)\in\Delta^\infty(\R)\cap\Delta^-(\R)$) such that $\mathcal
T^k(\la,\pi)\in\Delta^+(\R)$ (resp., $\mathcal
T^k(\la,\pi)\in\Delta^-(\R)$) for $k=0,\dots,n-1$, and $\mathcal
T^n(\la,\pi)\in\Delta^-(\R)$ (resp., $\mathcal
T^k(\la,\pi)\in\Delta^+(\R)$). Using (\ref{an}), one also can check
that if $w=(c,n,\pi)$, then
\begin{equation}
\label{Delta_w2} \Delta(w)=\begin{cases}
\Delta^\infty(\R)\cap\{(\la,\pi)|\,\la\in\Delta_\pi^+,\
\frac{A(w)^{-1}\la}{|A(w)^{-1}\la|}\in\Delta_{a^n\pi}^-\},\ \text{if}&c=a,\\
\Delta^\infty(\R)\cap\{(\la,\pi)|\,\la\in\Delta_\pi^-,\
\frac{A(w)^{-1}\la}{|A(w)^{-1}\la|}\in\Delta_{b^n\pi}^+\},\
\text{if}&c=b.
\end{cases}
\end{equation}
It is easy to see that all the sets $\Delta(w)$, $w\in\A$, are
non-empty and constitute a partition of $\Delta^\infty(\R)$. We use
this partition to construct a symbolic dynamics for $\mathcal G$. By
iterating this partition $n$ times under the action of the
transformation $\mathcal G$ we obtain a partition whose elements,
$\Delta(w)$, are determined by the words $w\in\wab$ of length $n$.
Namely, for a word $w=w_1\dots w_n\in\wab$, $w_i\in\mathcal A$, we
set
\begin{equation}
\label{Delta_w3} \Delta(w)=\bigcap_{i=0}^{n-1}\mathcal
G^{-i}\Delta(w_{i+1}).
\end{equation}
\begin{remark}
\label{Delta} From (\ref{te}), (\ref{ge}), and (\ref{renorm0})--
(\ref{Delta_w2}) it follows that if $(\la,\pi)\in\Delta(w)$ and
$(\la',\pi')=\cal G(\la,\pi)$, then
$$\la'=A^{-1}(w)\la/|A^{-1}(w)\la|,\ \ \la=A(w)\la'/|A(w)\la'|.$$
These formulas can be easily extended by induction to the case where
$w=w_1\dots w_n\in\wab$, $(\la,\pi)\in\Delta(w)$, and
$(\la',\pi')=\mathcal G(\la,\pi)$.
\end{remark}

The {\it coding map} $\Phi:\Delta^\infty(\R)\to\Omega_{\A, B}$ is
given by the formula
\begin{equation}
\label{cod} \Phi(\la,\pi)=(\omega_0,\omega_1,\dots)\ {\rm if} \
\gz^n(\la,\pi)\in\Delta(\omega_n)\ {\rm for\ }n=0,1,\dots.
\end{equation}

Consider conditions under which the coding map is invertible. For
$\q\in{\cal W}_{{\cal A},\,B}$ we denote by $\Omega_{\q}$ the set of
all sequences $\omega\in\omab$ starting from the word $\q$ and
containing infinitely many occurrences of $\q$. A key role will be
played below by the words $\q$ such that all entries of the
renormalization matrix $A(\q)$ are positive. For short we will refer
to these $\q$ as {\it positive words}. Observe that each word
containing a positive prefix is also positive.

The next two lemmas are due to Veech \cite{veech}.
\begin{lemma}
\label{recur} Let $\q\in{\cal W}_{{\cal A},\,B}$ be a positive word.
Then for every $\omega\in\Omega_{\q}$, there exists a unique
$(\la,\pi)\in\Delta^\infty(\R)$ such that $\Phi(\la,\pi)=\omega$.
\end{lemma}
\begin{lemma}
\label{uniq_erg} If the interval exchange map
$(\la,\pi):[0,1]\to[0,1]$ is uniquely ergodic, then there exists a
positive word $\q\in{\cal W}_{{\cal A},\,B}$ such that
$(\la,\pi)\in\Delta(\q)$.
\end{lemma}
\begin{remark}
\label{keane_pr} Veech \cite{veech} in fact observed that for the
existence of $\q$ mentioned in Lemma \ref{uniq_erg} it is sufficient
for $(\la,\pi)$ to satisfy Keane's infinite distinct orbit
condition.
\end{remark}

\subsection{Symbolic dynamics for the flow $\{P^t\}$}
\label{symb2}

We first construct a symbolic dynamics for the map $\mathcal F$
introduced in Subsection \ref{zipper}. For
$(\la,\pi,\delta)\in\mathcal Y_\infty^{\pm}(\R)$ we set
\begin{equation}
\label{tildePhi}
\tilde\Phi(\la,\pi,\delta)=(\dots,\omega_{-1},\omega_0,\omega_1,\dots),\
\ \omega_i\in\A,
\end{equation}
if $\mathcal F^n(\la,\pi,\delta)=(\la'_n,\pi'_n,\delta'_n)$ and
$(\la'_n,\pi'_n)\in\Delta(\omega_n)$, $n\in\mathbb Z$ (remind that
$\F$ is invertible on $\mathcal Y_\infty^{\pm}(\R)$). In parallel
with the coding map $\Phi$ (see Subsection \ref{symbol1}) we have
$$
\tilde\Phi(\Y_\infty^{\pm}(\R))\subset\Omega_{\A,\,B}^\mathbb Z.
$$
Moreover, from (\ref{ef}) and (\ref{cod}) it follows that
(\ref{tildePhi}) implies that
$$
\Phi(\la,\pi)=(\omega_0,\omega_1,\dots).
$$

For ${\bf q}\in\mathcal W_{{\mathcal A},B}$, $|{\bf q}|=l$,
similarly to the definition of $\Omega_{\q}$ in Subsection
\ref{symbol1}, we denote by $\omzq$ the set of all sequences
$\omega\in\omabz$ satisfying $\omega_0\dots \omega_{l-1}=\q$ and
admitting infinitely many occurrences, both in the past and in the
future, of the word $\q$.

Let
$$
\Y^{\pm}_{\q,\infty}(\R):={\tilde\Phi}^{-1}(\omzq),\ \ \ \ {\cal
V}_{\q}(\R):=\bigcup_{t\in{\mathbb R}}P^t\Y^{\pm}_{\q,\infty}(\R),
$$
and assume that $\Y_{\q,\infty}^{\pm}(\R)$ (and hence $\mathcal
V_\q(\R)$) is non-empty.

Let ${\F}_{\q}$ be the first return map of $\F$ to
$\Y^{\pm}_{\q,\infty}(\R)$, i.e., the map induced by $\F$ on
$\Y^{\pm}_{\q,\infty}(\R)$ (cf. Subsection \ref{induced_aut}).

By definition, $\{P^t\}|_{\mathcal V_{\bf q}(\R)}$, the restriction
of the flow $\{P^t\}$ to ${\cal V}_{\bf q}(\R)$, is Borel isomorphic
to a suspension flow $(\mathcal F_\q,\tau_\q)$ over the map
$\mathcal F_\q$. To describe the roof function $\tau_\q$ we take
$(\la,\pi)\in\Delta^\infty(\R)$, $(\lap,\pip)=\gz(\la,\pi)$, and
introduce, following Veech, the function
\begin{equation}
\label{deftau} \tau^1:(\la,\pi)\mapsto\log |A(\omega_0)\lap|,
\end{equation}
where $\omega_0$ is determined by the equation
$\Phi(\la,\pi)=(\omega_0,\omega_1,\dots)$, that is
$(\la,\pi)\in\Delta(\omega_0)$. (Recall that the norm of a vector
$v$ is given by $|v|=\sum_i |v_i|$.) Using (\ref{ef}), (\ref{tau2}),
it easy to check that if $x=(\la,\pi,\delta)\in\Y_\infty^{\pm}(\R)$
and $(\la,\pi)\in\Delta(\omega_0)$, where $\omega_0=(c,n,\pi)$, then
the first return time of $x$ to $\Y_\infty^{\pm}(\R)$ under the
action of the flow $\{P^t\}$ is just $\tau^1(\la,\pi)$.

Let now $(\la,\pi)\in\Phi^{-1}(\Omega_{\bf q})$,
$(\omega_0,\omega_1,\dots)=\Phi(\la,\pi)$, and let $s$ be the moment
of the second appearance of the word $\q$ in
$(\omega_0,\omega_1,\dots)$, i.e.,
\begin{equation}
\label{second}
s=s(\omega_0,\omega_1,\dots)=\min\{k>0|(\omega_k,\dots,\omega_{k+l-1})=\q\}.
\end{equation}
Denote
\begin{equation}
\label{tau4} \tau_\q^1(\la,\pi)=\tau^1(\la,\pi)+\tau^1(\mathcal
G(\la,\pi))+\dots+\tau^1(\mathcal G^{s-1}(\la,\pi)).
\end{equation}
If $x=(\la,\pi,\delta)\in\Y_{\bf q,\infty}^{\pm}(\R)$, then
$(\la,\pi)\in\Phi^{-1}(\Omega_{\bf q})$, and we have
\begin{equation}
\label{tau5} \tau_{\bf q}(x)=\tau_{\bf q}^1(\la,\pi).
\end{equation}

Denote by $\Psi_1$ the map from $\cal V_\q(\R)$ to the phase space
of $(\F_\q,\tau_q)$ that induces the above-mentioned isomorphism
between the flows $\{P^t\}|_{\cal V_\q(\R)}$ and $(\F_\q,\tau_q)$.

From now on we assume that $\q$ is a positive word. For such $\q$,
we construct a suspension flow $(\hat{\sigma}_\q, \hat\tau_{\q})$
closely related to $(\mathcal F_\q,\tau_\q)$. Denote, as before, by
$\sigma$ the one-step left shift on $\Omega_{{\mathcal
A},B}^{\mathbb Z}$ and, for
$\omega=(\dots,\omega_{-1},\omega_0,\omega_0,\dots)\in\Omega_{\bf
q}^{\mathbb Z}$, let
\begin{equation}
\label{tau6} \hat\sigma_{\bf q}(\omega):=\sigma^s,\ \ \hat\tau_{\bf
q}(\omega):=\tau_\q^1(\Phi^{-1}(\omega_0,\omega_1,\dots)),
\end{equation}
where $s$ is defined in (\ref{second}). (Observe that
$\omega\in\Omega_{\bf q}^{\mathbb Z}$ implies $(\omega_0,\omega_1, \
\dots)\in\Omega_{\bf q}$, hence by Lemma \ref{recur}
$\Phi^{-1}(\omega_0,\omega_1,\dots)$ is a uniquely defined point
from $\Delta^\infty(\R)$).

Proposition 6 in \cite{bufetov} states that if $\q\in{\cal W}_{{\cal
A},\,B}$ is positive, then for every $\omega\in\omzq$, there exists
at most one zippered rectangle corresponding to it; in other words,
the map ${\tilde \Phi}$ restricted to the set
$\Y_{\q,\infty}^\pm(\R):={\tilde\Phi}^{-1}(\omzq)$ is injective. It
induces, in a natural way, a measurable injective map $\Psi_2$ from
the phase space of the flow $(\sigma_\q,\tau_\q)$ to the phase space
of the flow $(\hat\sigma_\q,\hat\tau_\q)$ that sends the former flow
to the latter one restricted to some invariant set. Hence $(\mathcal
F_\q,\tau_\q)$ is embedded (in the sense of Section \ref{suspens})
into $(\hat\sigma_\q,\hat\tau_\q)$.

Introduce a new alphabet $\A_{\q}$; it will consist of all words
$w=(v_1\dots v_n)\in\wab$, $v_i\in\mathcal A$, $n>l$, such that
$(v_1,\dots,v_l)=\q$, $(v_{n-l+1},\dots v_n)=\q$, and no other
subword of $w$ coincides with $\q$. Since $|\A|=\infty$ and the
graph $\Gamma_{\A,B}$ is connected (see Subsection \ref{symbol1}),
we have $|\A_\q|=\infty$. By the Markov-Bernoulli reduction,
introduced in Subsection \ref{induced_aut}, there is a measurable
one-to-one map $\Psi_{M-B}:\Omega_\q^\mathbb Z\to(\mathcal
A_\q)^\mathbb Z$ that sends $\hat\sigma_\q$ to the one-step left
shift $\sigma_\q$ on $(\mathcal A_\q)^\mathbb Z$. Hence the flow
$(\hat{\sigma} _\q,\hat{\tau}_\q)$ is isomorphic to the suspension
flow $(\sigma_\q,f_\q)$, where
\begin{equation} \label{roof1}
f_\q(u):=\hat\tau_\q(\Psi_{M-B}^{-1}(u)),\ \ u\in(\mathcal
A_\q)^\mathbb Z.
\end{equation}
Denote by $\Psi_3$ the corresponding map from the phase space of
$(\hat{\sigma}_\q,\hat{\tau}_\q)$ to the phase space of
$(\sigma_\q,f_\q)$. Summing up, we can state the following.
\begin{lemma}
\label{embed} The mapping $\Psi:=\Psi_1\circ\Psi_2\circ\Psi_3$
yields an embedding of the flow $\{P^t\}|_{{\cal V}_\q(\R)}$ into
the flow $(\sigma_\q,f_\q)$.
\end{lemma}

Let us turn to the probability measure $\mu_\R$ on $\tilde{\cal
V}(\R)$ and the probability measure $\mu_\R^1$ induced by $\mu_\R$
on $\Y_\infty^\pm(\R)$ (see Subsection \ref{zipper}).

Observe that $\mu_{\R}$ assigns a positive mass to every Borel set
with nonempty interior. (We assume a natural topology on
$\tilde{\cal V}(\R)$ as well as on other spaces we encounter in this
paper.) From this fact, using Lemma \ref{recur} and the definition
of $\T$, $\cal G$, and $\F$, we easily derive that
$$
\mu^1_\R(\{x\in\Y_\infty^{\pm}(\R):x=(\la,\pi,\delta),\,
(\la,\pi)\in\Delta(\q)\})>0.
$$
Since $\mu_{\R}$ is $\{P^t\}$-ergodic, we have
\begin{equation}
\label{full_meas} \mu_{\R}({\cal V}_{\q}(\R))=1,\ \
\mu^1_\R\left(\bigcup_{n\in\mathbb
Z}\F^n\Y_{\q,\infty}^\pm(\R)\right)=1.
\end{equation}
By normalizing the restriction of $\mu_\R^1$ to
$\Y_{\q,\infty}^\pm(\R)$ we obtain a probability measure
$\bar\mu^1_{\R,\q} $.

Let $\psi$ be the natural projection of $\omabz$ on $\omab$. Recall
that in Subsection \ref{zipper} we introduced the natural projection
$\tilde\psi$ of $\Y_\infty^\pm(\R)$ on $\Delta^\infty(\R)$. It is
clear that $\psi(\tilde\Phi(x))=\Phi(\tilde\psi(x))$ for every
$x\in\bigcup_{n\in\mathbb Z}\F^n \Y_{\q,\infty}^\pm(\R)$. From this
fact combined with (\ref{full_meas}) and Lemma \ref{meas_proj} we
obtain
\begin{equation}
\label{proj_nu} \psi_*(\tilde\Phi_*\mu_\R^1)=\Phi_*\nu.
\end{equation}

Using the positivity of $\nu$ on all open sets, we come to the
following assertion.
\begin{lemma} \label{pos_mes} The measure
$\Psi_*\mu_\R$ is positive on every open set in the phase space of
the flow $(\sigma_\q,f_\q)$.
\end{lemma}

\subsection{Properties of the roof function}
\label{roof}

Recall that a word $w'=w_1\dots w_l\in\W_{\A,B}$ is said to be a
simple prefix of a word $w=w_1,\dots,w_l,\dots,w_n$ if $w_1\dots
w_{n-k+1}=w_k\dots w_n$ implies that either $k=1$, or $k>l$ (see
Subsection \ref{over_markov}).
\begin{lemma}
\label{hoelder} Let $\q\in\W_{\A,B}$ be a word that has a simple
positive prefix. Then the function $f_\q$ introduced in
(\ref{roof1}) depends only on the future, is bounded away from zero,
and is H\" older continuous in the following sense: there exist
positive constants $C_{\q}$, $\alpha_{\q}$ (depending only on $\q$)
such that if $u=(\dots,u_{-1},u_0,u_1\dots)\in\A_{\q}^{{\mathbb Z}}$
and ${\tilde u}=(\dots,{\tilde u}_{-1},{\tilde u}_0,{\tilde
u}_1\dots)\in\A_{\q}^{{\mathbb Z}}$ satisfy $u_i={\tilde u}_i$ for
$|i|\leq n$, then
$$
|f_{\q}(u)-f_{\q}(\tilde u)|\leq C_{\q}\exp(-\alpha_{\q}n).
$$
In particular, the function $f_{\q}$ has summable variations.
\end{lemma}
\begin{proof} That the function $f_\q$ depends only on the future
follows readily from its definition. Let us prove that it is bounded
away from zero. By (\ref{roof1}), (\ref{tau6}) this property of
$f_\q$ would follow from the same property of the function
$\tau_\q^1$ defined on $\Phi^{-1}(\Omega_\q)$. But from
(\ref{deftau}), (\ref{tau4}) and Remark \ref{Delta} one readily
derives that, for $(\omega_0,\omega_1,\dots)\in\Omega_\q$,
\begin{equation}
\label{tau4a} \tau_\q^1(\la,\pi)=\log|A(\omega_0)\dots
A(\omega_{s-1})\la'|,\ \
(\la,\pi)=\Phi^{-1}(\omega_0,\omega_1,\dots),
\end{equation}
where $s$ is defined in (\ref{second}) and $(\la',\pi')=\mathcal
G(\la,\pi)=\Phi^{-1}(\omega_s,\omega_{s+1},\dots)$. If $\p$ is a
simple prefix of $\q$, then $s\ge|\p|$, so that $\p$ is a prefix of
the word $(\omega_0,\dots,\omega_{s-1})$ and hence this word is
positive. Thus all entries of the $m\times m$ matrix
$$A(\omega_0,\dots,\omega_{s-1})=A(\omega_0)\dots
A(\omega_{s-1})
$$
are positive integers, while $\la'$ is a positive vector with
$|\la'|=1$. Therefore, for every $(\la,\pi)\in\Phi^{-1}(\Omega_\q)$,
$$\tau_\q^1(\la,\pi)=\log|A(\omega_0,\dots,\omega_{s-1})\la'|\ge\log m>0.$$

We now turn to the H\"older continuity of $f_\q$. Let $u,\tilde u$
be as in the statement and let $\omega=\Psi^{-1}_{M-B}(u)$,
$\tilde\omega=\Psi^{-1}_{M-B}(\tilde u)$. By definition,
$\omega,\tilde\omega\in\Omega_\q^{\mathbb Z}$, hence
$(\omega_0,\omega_1,\dots),\,(\tilde\omega_0,\tilde\omega_1,\dots)\in\Omega_
\q$. By (\ref{tau5})--(\ref{roof1})
\begin{equation}
\label{tau4b}
f_\q(u)=\tau_\q^1(\la,\pi)=\log|A(\omega_0,\dots,\omega_{s-1})\la'|,
\end{equation}
\begin{equation}
\label{tau4c}
 f_\q(\tilde u)=\tau_\q^1(\tilde\la,\tilde\pi)=
\log|A(\tilde\omega_0,\dots,\tilde\omega_{s-1})\tilde\la'|,
\end{equation}
where
$$(\la,\pi)=\Phi^{-1}(\omega_0,\omega_1,\dots),\ \
(\tilde\la,\tilde\pi)=\Phi^{-1}(\tilde\omega_0,\tilde\omega_1,\dots),
$$
$$
(\la',\pi')=\mathcal G^s(\la,\pi),\ \
(\tilde\la',\tilde\pi')=\mathcal G^s(\tilde\la,\tilde\pi).
$$
Since $\p$ is a simple prefix of $\q$, one can find $k\ge nl-1$ such
that $\omega_i=\tilde\omega_i$ for $i=0,1,\dots,k$. Then, by Remark
\ref{Delta}, for some vectors $\la'',\tilde\la''\in\Delta_{m-1}$,
\begin{equation}
\label{lambda} \la=\frac{A(\omega_0,\dots\omega_{nl-1})\la''}
{|A(\omega_0,\dots\omega_{nl-1})\la''|},\ \ \
\tilde\la=\frac{A(\omega_0,\dots\omega_{nl-1})\tilde\la''}
{|A(\omega_0,\dots\omega_{nl-1})\tilde\la''|},
\end{equation}

Introduce the Hilbert metric $d_H$ on $\Delta_{m-1}$ by
\begin{equation}
\label{hilbert} d_H(\la^{(1)},\la^{(2)})=\log\left(\max_{1\le i\le
m}\frac{\la_i^{(1)}}{\la_i^{(2)}}/\min_{1\le i\le
m}\frac{\la_i^{(1)}}{\la_i^{(2)}}\right),\
\la^{(j)}=(\la^{(j)}_1,\dots,\la^{(j)}_m), j=1,2.
\end{equation}
It is known that if an $m\times m$ non-negative matrix $A=(a_{ij})$
is such that $\sum_{j=1}^ma_{ij}>0$ for all $i$, then the mapping
$T_A:\Delta_{m-1}\to\Delta_{m-1}$ defined by $T_A\la=A\la/|A\la|$
does not increase the $d_H$-distance between points, while if
$a_{ij}>0$ for all $i,j$, then $T_A$ is a uniform contraction (see,
for example, \cite{viana}).

By definition, the word $\omega_0\dots\omega_{nl-1}$ is a
concatenation, namely, $\omega_0\dots\omega_{nl-1}=\p w_1\dots\p
w_n$, where the word $w_i$, $1\le i\le n$, is such that the sum of
the entries in each row of the matrix $A(w_i)$ is positive. It
follows (see (\ref{lambda})) that
\begin{equation}
\label{hilbert1} d_H(\la,\tilde\la)\le
d_H(T_{A(\q)}^n\la'',T_{A(\q)}^n\tilde\la'')\le C_1\alpha^n,
\end{equation}
where $C_1\in\mathbb R_+$ and $\alpha\in(0,1)$ depend only on $\q$.
(We have used the fact that $T_{A(\q)}$ takes $\Delta_{m-1}$ to a
set of finite $d_H$-diameter).

Denote $A_1=A(\omega_0,\dots,\omega_{s-1})$. Since all entries of
the matrix $A_1$ are positive, we have (see (\ref{hilbert}))
$$d_H\left(\frac{A_1\la}{|A_1\la|},
\frac{A_1\tilde\la}{|A_1\tilde\la|}\right)=
\log\left(\max_i\frac{(A_1\la)_i}{(A_1\tilde\la)_i}/
\min_i\frac{(A_1\la)_i}{(A_1\tilde\la)_i}\right)\le
C_2d_H(\la,\tilde\la),
$$
where $(A_1\la)_i$ (resp., $(A_1\tilde\la)_i$) is the $i$th entry of
the vector $A_1\la$ (resp., $A_1\tilde\la$) and $C_2$ is determined
by $\q$. Hence, by (\ref{tau4b}) and (\ref{hilbert1}),
$$|f_\q(u)-f_\q(\tilde u)|=|\log(|A_1\la|/|A_1\tilde\la|)|\le
C_2d_H(\la,\tilde\la)\le C_2C_1\alpha^n,
$$
so it remains to set $C_\q=C_2C_1$ and $\alpha_\q=-\log\alpha$.
\end{proof}

\subsection{Transition probabilities and the uniform expansion property.}
\label{transit}

By Theorem \ref{zorichthm}, the map $\gz$ on $\Delta(\R)$ preserves
an absolutely continuous ergodic probability measure, which was
denoted by $\nu$.

Consider a word $w=w_1\dots w_k\in{\cal W}_{\A,B}$, where
$w_i=(c_i,n_i,\pi_i)\in\mathcal A$, $1\le i\le k$. We say that $w$
is {\it compatible} with a point $(\la,\pi)\in\Delta^\infty(\R)$ (or
$(\la,\pi)$ is compatible with $w$) if
$$\text{either }\la\in\Delta_{\pi}^{-},\ c_k=a,\
a^{n_k}\pi_k=\pi,\text{\ \ or\ \ }\la\in\Delta_{\pi}^{+},\ c_k=b,\
b^{n_k}\pi_k=\pi.
$$
Assuming that $w$ is compatible with $(\la,\pi)$, we set
\begin{equation}
\label{te_w}
t_w(\la,\pi)=\left(\frac{A(w)\la}{|A(w)\la|},\pi_1\right).
\end{equation}
From the definition of $\cal G$ (see (\ref{ge})) it follows that
\begin{equation}
\label{te_w1} \gz^{-n}(\la,\pi)=\{t_w(\la,\pi):|w|\ {\rm and}\ w\
{\rm is\ compatible \ with\ } \ (\la,\pi)\}.
\end{equation}
Note that the set $\gz^{-n}(\la,\pi)$ is infinite.

In Subsection \ref{symbol1} we introduced, for each word $w\in\wab$,
the set $\Delta(w)\subset\Delta^\infty(\R)$ (see (\ref{Delta_w3})).
One can readily check (see Remark \ref{Delta}) that
\begin{equation}
\label{Delta_w4} \Delta(w)=\{t_w(\la,\pi): \
(\la,\pi)\in\Delta^\infty(\R)\ {\rm\ is\ compatible\ with}\ w\}.
\end{equation}

For every $n\in\mathbb N$, we have the $\nu$-measurable partition
$\mathcal G^{-n}\varepsilon$ of $\Delta^\infty(\R)$, where
$\varepsilon$ is the partition into separate points. Each element of
$\mathcal G^{-n}\varepsilon$ is $\mathcal G^{-n}(\la,\pi)$ for some
$(\la,\pi)\in\Delta^\infty(\R)$, its points correspond to the words
$w\in\mathcal W_{\mathcal A,B}$ of length $n$ compatible with
$(\la,\pi)$ and have the form $t_w(\la,\pi)$ (see (\ref{te_w}),
(\ref{te_w1})). We denote by $\nu(w|(\la,\pi))$ the conditional
measure (determined by $\nu$) of the point corresponding to $w$,
given the element $\mathcal G^{-n}(\la,\pi)$ of the partition
$\mathcal G^{-n}\varepsilon$.

In Section 3.5 of \cite{bufetov} it is proved that if $w$ is
compatible with $(\la,\pi)$, then
\begin{equation}
\label{trpr}
\nu(w|(\la,\pi))=\frac{\rho(t_w(\la,\pi))}{\rho(\la,\pi)|A(w)\la|^m}.
\end{equation}

Now consider the set $\Delta(\q)$ corresponding to a word
$\q\in\wab$ (see (\ref{Delta_w3})). Every point from $\Delta(\q)$ is
of the form $(\la,\pi_{\q})$, where $\pi_{\q}$ is a fixed
permutation and $\la$ belongs to a set $\Delta^{\prime}(\q)\subset
\Delta_{m-1}$. Denote by ${\rm d}(\q)$ the diameter of
$\Delta^{\prime}(\q)$ with respect to the Hilbert metric on
$\Delta_{m-1}$ introduced in (\ref{hilbert}).
\begin{proposition}
\label{unifexp} There are positive constants $\beta_1$ and $\beta_2$
(depending only on $\R$) such that for every positive word
$\p'\in\wab$ with ${\rm d}(\p')\le\beta_1$, the following holds. Let
$\p$ be a word that has $\p'$ as a simple prefix, and let a word
${\bf r}\in{\cal W}_{\A,B}$ start and end with $\p$ and contain no
other occurrences of $\p$. Then for any $(\la,\pi)\in\Delta({\bf
r})\cap\Phi^{-1}(\Omega_\p)$, we have
$$
\left|\frac{\nu(\Delta({\bf r}))\exp[m\tau_{\p}^1(\la,\pi)]}
{\nu(\Delta(\p))}-1\right|\leq\beta_2{\rm d}(\p'),
$$
where $\tau_{\p}^1$ is defined by (\ref{tau4a}).
\end{proposition}
\begin{proof}
By assumption, the word ${\bf r}$ has the form ${\bf r}=\p'
u\p=\p{\tilde u}$ for some $u, {\tilde u}\in\wab$. Hence
\begin{equation}
\label{qwq} \nu(\Delta({\bf r}))=\int_{\Delta(\p)}\nu(\p'
u|(\la,\pi))d\nu(\la,\pi).
\end{equation}

From (\ref{trpr}) we obtain
\begin{equation}
\label{invdens} \nu(\p' u|(\la,\pi))=\frac{\rho(t_{\p' u}(\la,\pi))}
{\rho(\la,\pi)}\cdot\frac 1{|A(\p' u)\la|^m}.
\end{equation}

By (\ref{hilbert}), taking into account that $|\la|=1$ when
$\la\in\Delta_{m-1}$, for any two points,
$(\la,\pi_{\p'}),(\tilde\la,\pi_{\p'})\in\Delta(\p')$, and for
$i=1,\dots,m$, we have
\begin{equation}
\label{diam1} e^{-{\rm d}(\p')}\leq\tilde\la_i/\la_i\leq e^{{\rm
d}(\p')}.
\end{equation}

Let us estimate the first ratio on the right-hand side of
(\ref{invdens}). From the fact that $\p' u\p'\in\wab$ it follows
that if $(\la,\pi_{\p'})\in\Delta(\p')$, then $(\la,\pi_{\p'})$ is
compatible with $\p' u$. Hence $t_{\p' u}\in\Delta(\p'
u)\subset\Delta(\p')$ (see (\ref{Delta_w4})). Denote
$$
(\hat\la,\hat\pi):=t_{\p' u}(\la,\pi_{\p'}),\ \
\hat\la=(\hat\la_1,\dots,\hat\la_m),
$$
and observe that $\hat\pi=\pi_{\p'}$. By (\ref{diam1})
$$
e^{-{\rm d}(\p')}\la_i\le\hat\la_i\le e^{{\rm d}(\p')}\la_i,\ \
i=1,\dots,m,
$$
while by Theorem \ref{zorichthm}
$\rho(\hat\la,\hat\pi)=P_{\hat\pi}(\hat\la)/Q_{\hat\pi}(\hat\la)$,
where $P_{\hat\pi}$ and $Q_{\hat\pi}$ are homogeneous polynomials
with non-negative coefficients. Therefore,
$$
P_{\hat\pi}(\hat\la)\le P_{\hat\pi}(e^{{\rm d}(\p')}\la)\le
e^{\gamma_1{\rm d}(\p')}P_{\hat\pi}(\hat\la),
$$
$$
Q_{\hat\pi}(\hat\la)\ge Q_{\hat\pi}(e^{-{\rm d}(\p')}\la)\ge
e^{-\gamma_2{\rm d}(\p')}Q_{\hat\pi}(\hat\la),
$$
where $\gamma_1$ and $\gamma_2$ are determined by $\rho$ (and,
eventually, by $\R$). From this we immediately obtain
\begin{equation}
\label{rat} \exp[-(\gamma_1+\gamma_2){\rm
d}(\p')]\le\frac{\rho(t_{\p' u}(\la,\pi))}{\rho(\la,\pi)}\le
\exp[(\gamma_1+\gamma_2){\rm d}(\p')].
\end{equation}

Using (\ref{diam1}) and the positivity of all entries of the matrix
$A(\p' u)$, for any $(\la,\pi_{\p'}),
(\tilde\la,\pi_{\p'})\in\Delta(\p')$, we have
\begin{equation}
\label{ctt} \exp[-{\rm d}(\p')]\le \frac{|A(\p' u)\tilde\la|}{|A(\p'
u)\la|}\le\exp[{\rm d}(\p')].
\end{equation}

Now fix an arbitrary point $(\la^0,\pi^0)\in\Delta({\bf
r})\cap\Phi^{-1}(\Omega_\p)$. By Lemma \ref{recur}
$(\la^0,\pi^0)=\Phi^{-1}(\omega_0,\omega_1,\dots)$ for some
$(\omega_0,\omega_1,\dots)\in\Omega_\p$. Let $s$ be defined by
(\ref{second}) and $(\la',\pi'):=\mathcal G^s(\la^0,\pi^0)$.
Clearly, $s=|\p' u|$, hence
$(\la',\pi')\in\Delta(\p)\subset\Delta(\p')$, and by (\ref{qwq}),
(\ref{invdens}), (\ref{rat}), (\ref{ctt}),
\begin{align}
\label{ratio1} \frac{\nu(\Delta({\bf r}))|A(\p'
u)\la'|^m}{\nu(\Delta(\p))}&=
\frac{1}{\nu(\Delta(\p))}\underset{\Delta(\p)}{\int}
{\frac{\rho(t_{\p' u}(\la,\pi))}
{\rho(\la,\pi)}\cdot\frac{|A(\p' u)\la'|^m}{|A(\p' u)\la|^m}\,d\nu(\la,\pi)}\notag\\
&\le\exp[(\gamma_1+\gamma_2+1){\rm d}(\p')].
\end{align}
Similarly,
\begin{equation}
\label{ratio2} \frac{\nu(\Delta({\bf r}))|A(\p'
u)\la'|^m}{\nu(\Delta(\p))}\ge \exp[-(\gamma_1+\gamma_2+1){\rm
d}(\p')].
\end{equation}
From (\ref{ratio1}), (\ref{ratio2}) we obtain
$$
\left|\frac{\nu(\Delta({\bf r}))|A(\p'
u)\la'|^m}{\nu(\Delta(\p))}-1\right|\le
\exp[(\gamma_1+\gamma_2+1){\rm d}(\p)]\le 2(\gamma_1+\gamma_2+1){\rm
d}(\p'),
$$
where the last inequality holds when $(\gamma_1+\gamma_2+1){\rm
d}(\p')\le\log 2$. It remains to recall that by (\ref{tau4a})
$|A(\p' u)\la'|=\exp(\tau_\p^1(\la,\pi))$.
\end{proof}

\section{Zippered rectangles and Abelian differentials. Completion of the
proof of Theorem \ref{maxentropy}} \label{abel}

Fix a connected component $\HH$ of the space $\modk$ (see Section
\ref{introd}). To this component there corresponds a unique Rauzy
class $\R$ in such a way that the following is true \cite{veech,
KZ}.

\begin{theorem}[\rm Veech]
\label{zipmodule} There exists a finite-to-one measurable map
$\pi_{\R}:\tilde{\cal V}(\R)\to\HH$ such that $\pi_{\R}\circ
P^t=g_t\circ \pi_{\R}$ and $(\pi_{\R})_*\mu_{\R}=\mu_{\kappa}$ for
all $t\in\mathbb R$.
\end{theorem}
\noindent(Recall that the set $\tilde{\cal V}(\R)$ is defined in
Subsection \ref{zipper}).
\begin{corollary}
\label{correspondence-measures} 1. If $\eta$ is a
$\{g_t\}$-invariant ergodic probability measure on $\HH$, then there
exists a $\{P^t\}$-invariant measure $\tilde\eta$ on
$\tilde{\mathcal V}(\R)$ with $(\pi_{\R})_*\tilde\eta=\eta$.

2. If $\tilde\eta$ is a $\{P^t\}$-invariant probability measure on
$\tilde{\mathcal V}(\R)$ such that the $\{g_t\}$-invariant measure
$(\pi_{\R})_*\tilde\eta$ is ergodic, then
\begin{equation}
\label{ent}
h_{\tilde\eta}(\{P^t\})=h_{(\pi_{\R})_*\tilde\eta}(\{g_t\}).
\end{equation}
\end{corollary}
\begin{proof}
1. Let $\eta$ be an ergodic $\{g_t\}$-invariant probability measure
on $\HH$. By ergodicity, there is a set $\HH'\subset\HH$ such that
$\eta({\HH'})=1$ and the cardinality of the preimage
$\pi_{\R}^{-1}(p)$ does not depend on $p\in\HH'$. The sets
$\pi_{\R}^{-1}(p)$ form a measurable partition of $\tilde{\mathcal
V}(\R)$. By assigning equal weights (conditional measures) to all
points in $\pi_{\R}^{-1}(p)$, $p\in\HH$, we obtain a
$\{P^t\}$-invariant probability measure $\tilde\eta$ on
$\tilde{\mathcal V}(\R)$ such that $(\pi_{\R})_*\tilde\eta=\eta$.

2. Let $\tilde\eta$ be a $\{P^t\}$-invariant probability measure on
$\tilde{\mathcal V}(\R)$ and $\eta=(\pi_{\R})_*\tilde\eta$. Assume
that $\eta$ is ergodic with respect to the flow $\{g_t\}$, hence it
is ergodic with respect to the automorphism $g_{t_0}$ for some
$t_0>0$.

Denote by $\frak c$ the canonical partition for $P^{t_0}$ and
$\tilde\eta$, i.e., the partition of $\tilde{\mathcal V}(\R)$
corresponding to the decomposition of $P^{t_0}$ into ergodic
components with respect to $\tilde\eta$. For an element $C$ of
$\frak c$, denote by $\tilde\eta_C$ the conditional measure induced
by $\tilde\eta$ on $C$, and consider $\tilde\eta_C$ as a measure on
$\tilde{\mathcal V}(\R)$. Let us note that for $\tilde\eta$-almost
all $C$, the measure $\tilde\eta_C$ exists and is ergodic with
respect to $P^{t_0}$.

Take a measurable set $M$ comprised of preimages $\pi_\R^{-1}(p)$.
From the Birkhoff ergodic theorem, applied to $P^{t_0}$,
$\tilde\eta$, $M$ and to $P^{t_0}$, $\tilde\eta_C$, $M$, it follows
that $\frak c$ and the partition of $\tilde{\mathcal V}(\R)$ into
the preimages $\pi_{\R}^{-1}(p)$, $p\in\mathcal H$, are independent
with respect to the measure $\tilde\eta$. Therefore
$(\pi_{\R})_*\tilde\eta_C=\eta$ for $\tilde\eta$-almost all elements
$C\in\frak c$.

As before, the cardinality of $\pi_{\R}^{-1}(p)$ is a constant, say
$k$, on a set $\mathcal H'\subset\mathcal H$ with $\eta(\mathcal
H')=1$. Moreover, if $\tilde\eta_C$ is $P^{t_0}$-ergodic, then the
conditional measure induced by $\tilde\eta_C$ on $\pi_{\R}^{-1}(p)$
for $\eta$-almost all $p\in\mathcal H'$ is uniform. Therefore, the
measure space $(\tilde{\mathcal V}(\R),\tilde\eta_C)$ is isomorphic
to the direct product of $(\mathcal H,\eta)$ and a set consisting of
$k$ points of mass $1/k$ each. In this representation, $P^{t_0}$
becomes a skew product with base $(\mathcal H,\eta,g_{t_0})$, and
the Abramov-Rokhlin formula \cite{AR} for the entropy of a skew
product implies that $h_{\tilde\eta_C}(P^{t_0})=h_\eta(g_{t_0})$.
Thus the entropy of $\tilde\eta$-almost every ergodic component of
$P^{t_0}$ equals $h_\eta(g_{t_0})$, which implies that
$h_{\tilde\eta}(P^{t_0})=h_\eta(g_{t_0})$. Since the entropy of
automorphisms forming a measurable flow $\{S_t\}$ with respect to an
$\{S_t\}$-invariant probability measure $\mu$ satisfies the equation
$h_\mu(\{S_t\})=|t|h_\mu(S_1)$, we come to (\ref{ent}).
\end{proof}

From Theorem \ref{zipmodule} and Corollary
\ref{correspondence-measures}, taking into account the ergodicity of
the measure $\mu_\kappa$ with respect to $\{g_t\}$ (see Section
\ref{introd}), we readily obtain the following.
\begin{corollary}
\label{mu_max} To prove that $\mu_{\kappa}$ is a unique measure with
maximal entropy for the flow $\{g_t\}$ it suffices to show that
$\mu_{\R}$ is a unique measure with maximal entropy for the flow
$\{P^t\}$.
\end{corollary}

Recall that the Rauzy class $\R$ we deal with consists of
permutations on $m$ symbols. With this in mind we derive from
(\ref{entropy}) and Corollary \ref{correspondence-measures} that
\begin{equation}
\label{entropy1} h_{\mu_{\R}}({\{P^t\}})=2g-1+r=m,
\end{equation}

We call a point $x\in\tilde{\mathcal V}(\R)$ {\it infinitely
renormalizable} if its trajectory $\{P^tx,\,t\in\mathbb R\}$
intersects the transversal $\Y^{\pm}(\R)$ infinitely many times both
for $t>0$ and for $t<0$. The set of infinitely renormalizable points
was denoted in Section \ref{zipper} by $\Y_\infty^\pm(\R)$.

The following proposition is in essence contained in \cite{veech}
and \cite{masur}.
\begin{proposition}
\label{set_V} There exists a Borel measurable set $V\subset\mathcal
H$ such that

(i) $\mu(V)=1$ for every ergodic $\{g_t\}$-invariant probability
measure $\mu$ on $\mathcal H$;

(ii) $\pi_\R^{-1}(p)\cap{\Y}_\infty^\pm\ne\emptyset$ for each $p\in
V$.
\end{proposition}
\begin{proof}
For a compact set $K\subset\mathcal H$, denote by $K^\pm$ the set of
points $p\in\mathcal H$ for which there exist $t_n\to+\infty$ and
$s_n\to-\infty$ such that $g_{t_n}p\in K$, $g_{s_n}y\in K$ for
$n=1,2,\dots$.

Take an increasing sequence of compact sets $K_n$ such that
$\underset{n}{\bigcup}K_n=\mathcal H$, and let
$V=\underset{n}{\bigcup}K_n^\pm$. The set $V$ is obviously Borel
measurable. By definition, for every probability measure $\mu$ on
$\mathcal H$, there exists $n_0$ with $\mu(K_{n_0})>0$. If, in
addition, $\mu$ is $\{g_t\}$-invariant and ergodic, then
$\mu(K_{n_0}^\pm)=1$ and hence $\mu(V)=1$.

Let $(\sigma,\omega)$ belong to an equivalence class $p\in V$ (see
Section \ref{introd}). Then, by Masur's theorem \cite{masur}, the
foliations corresponding to $\Re(\omega)$ and $\Im(\omega)$ are both
uniquely ergodic. This implies, in particular, the existence of an
infinitely renormalizable zippered rectangle in $\pi_\R^{-1}p$ (see,
for instance, \cite{zorich1}), which is all that we had to prove.
\end{proof}
\begin{lemma}
\label{positive} Let $\tilde\eta$ be an ergodic $\{P^t\}$-invariant
probability measure on $\tilde{\mathcal V}(\R)$. Then there exists a
positive word $\p\in\wab$ such that $\tilde\eta({\cal
V}_{\p}(\R))=1$. Moreover, if $h_{\tilde\eta}(\{P^t\})>0$, then for
any ${\p}$ such that $\tilde\eta({\cal V}_{\p}(\R))=1$, there also
exists a word $\q\in\wab$ such that $\tilde\eta({\cal
V}_{\q}(\R))=1$ and $\p$ is a simple prefix of $\q$.
\end{lemma}
\begin{proof}
Since the measure $\tilde\eta$ is ergodic with respect to $\{P^t\}$,
its projection $\eta:=(\pi_\R)_*\tilde\eta$ is ergodic with respect
to $\{g_t\}$. Let $V$ be as in Proposition \ref{set_V} and
$W=\pi_{\R}^{-1}(V)$. By Proposition \ref{set_V} $\eta(V)=1$, which
implies that $\tilde\eta(W)=1$.

For every point $x=(\la,\pi,\delta)\in W$, the interval exchange
$(\la,\pi)$ is uniquely ergodic. Now Lemma \ref{uniq_erg} yields a
positive word $\q_x\in\wab$ such that $(\la,\pi)\in\Delta(\q_x)$.
Since $\eta(W)=1$, while $\wab$ is countable, there exists
$\q\in\wab$ such that $\tilde\eta(\{x:\q_x=\q\})>0$. For this $\q$
we have $\tilde\eta({\mathcal Y}_{\q,\infty}^\pm)(\R)>0$ and, since
$\tilde\eta$ is ergodic, $\tilde\eta(\mathcal V_\q(\R))=1$.

Let us now consider the set $\mathcal W(\p)\subset\wab$ of all words
that have $\p$ as a prefix. If $\q\in\mathcal W(\p)$, and $\p$ is
not a simple prefix of $\q$, then $\q$ is a concatenation:
$\q=\q'\q'\dots\q'\q''$, where $\q'$ is a prefix of $\p$, and $\q'$
is either a prefix of $\q'$ or empty. In this situation, either
$\tilde\eta(\mathcal V_\q(\R))=0$ or $\tilde\eta(\mathcal
V_\q(\R))>0$ and $\tilde\eta$ is concentrated on periodic points of
the flow $\{P^t\}$. But this can not be the case for all
$\q\in\mathcal W(\p)$, since $\tilde\eta(\mathcal V_\q(\R))=1$ and
$h_{\tilde\eta}(\{P^t\})>0$.
\end{proof}

The following statement will be also used below.
\begin{lemma}
\label{longsimple} Let $\Gamma=(V,E)$ be a directed graph with
$|V|=\infty$ and $W(\Gamma)$ be the corresponding family of words
(see Subsection \ref{over_markov}). Then for each $w\in W(\Gamma)$
and each $n\in\mathbb N$, there exists a word $w'\in W(\Gamma)$ of
the form $w'=ww_1ww_2w,\dots w_{n-1}w$, $w_i\in W(\Gamma)$,
containing $n$ disjoint subwords equal to $w$ and such that the word
$w'':=ww_1ww_2w,\dots w_{n-1}$ is simple.
\end{lemma}
\begin{proof}
Denote the first and last letters of $w$ by $v^-$ and $v^+$
respectively, and construct by induction a sequence of letters
$v_1$, $v_2$,\dots as follows. For $v_1$ we take an arbitrary letter
that is not contained in $w$. If $v_1,\dots,v_k$, $1\le k\le n-1$,
are already chosen, we denote by $w_k^+$ (respectively, $w_k^-$) the
shortest word from $v^+$ to $v_k$ (resp., from $v_k$ to $v^-$). (If
there are several words of the same, minimal, length, we take any of
them.) Then take for $v_{k+1}$ an arbitrary vertex except for those
contained in at least one of the words $w,
w_1^+,w_1^-,\dots,w_k^+,w_k^-$.

Consider the sequence of words
$$
w,w_1^+,w_1^-,w,w_2^+,w_2^-,w,\dots,w,w_n^+,w_{n-1}^-,w.
$$
Delete the last letter from each of these words but the last $w$,
and denote by $w'$ the concatenation of the words thus obtained. It
is easy to check that $w'$ possesses the required properties. (To
define $w_k$ one should remove the first and last letters from
$w_k^+$, the last letter from $w_k^-$, and then take the
concatenation of the two words obtained.)
\end{proof}
\smallskip

{\bf End of the proof of Theorem \ref{maxentropy}}. By Corollary
\ref{mu_max}, it suffices to prove that $\mu_\R$ is a unique measure
with maximal entropy for the flow $\{P^t\}$.

Let $\mu$ be a measure on $\tilde V(\R)$ with $h_\mu(\{P^t\})\ge
h_{\mu_\R}(\{P^t\})$. Without loss of generality we can assume that
$\mu$ is ergodic (otherwise one could pass to an ergodic component).
By Lemma \ref{positive}, there exists a word $\q\in\mathcal
W_{\mathcal A,B}$ that has a simple positive prefix and is such that
$\mu(\mathcal V_\q(\R))=1$. Remind that $\mu_\R(\mathcal
V_\q(\R))=1$ as well. By Lemma \ref{embed}, the flow
$\{P^t\}|_{\mathcal V_\q(\R)}$ is embedded in the suspension flow
$(\sigma_\q,f_\q)$ via a mapping $\Psi$. The roof function $f_\q$
clearly depends only on the future (see Subsection
\ref{over_markov}). Moreover, by Lemma \ref{hoelder} $f_\q$ is
bounded away from zero and has summable variations.

By Lemma \ref{pos_mes} the measure $\Psi_*\mu_\R$ is positive on all
non-empty open subsets of the phase space of the flow
$(\sigma_\q,f_\q)$. This measure induces, in a canonical way, a
probability measure $\bar\mu_\q$ on $(\mathcal A_\q)^{\mathbb Z}$,
the base of the suspension flow $(\sigma_\q,f_\q)$. It follows that
$\bar\mu_\q$ is positive on all cylinders in $(\mathcal
A_\q)^{\mathbb Z}$.

Let us prove that ${\bar\mu}_\q$ satisfies the other conditions
imposed on the measure $\bar\mu$ in Theorem \ref{main}. (We apply
this theorem to the complete graph with vertex set $\mathcal A_\q$;
in this situation every sequence of letters (i.e., vertices of the
graph) is a word.)

Apply Lemma \ref{longsimple} to the connected graph $\Gamma_{{\cal
A},B}$ (see Subsection \ref{symbol1}) and to the word $\q$ taken as
$w$. From this lemma we obtain, for each $n\in\mathbb N$, a word of
the form $\q q_1\q\dots\q q_n\q\in\wab$ whose prefix $\q
q_1\q\dots\q q_n$ is a simple word.

Each word $\ba_i:=\q q_i\q$, $i=1,\dots,n$, is a letter in the
alphabet $\mathcal A_\q$ introduced in Subsection \ref{symb2}, and
$\ba:=(\ba_1,\dots,\ba_n)$ is clearly a simple word for all $n$.
Denote by $\hat a$ an arbitrary word in the alphabet $\A_\q$ that
does not contain $\ba$ as a subword, and consider the cylinders
$C_\ba$ and $C_{\ba\hat a\ba}\subset(\A_\q)^\mathbb Z$. Using the
definition of the measures and maps that appear below, we have
\begin{equation}
\label{ratio3} \frac{\bar\mu_\q(C_{\ba\hat
a\ba})}{\bar\mu_\q(C_\ba)}=\frac{(\tilde\Phi_*\mu_{\R,\q}^1)(\Psi_{M-B}^{-1}C_{\ba\hat
a\ba})}{(\tilde\Phi_*\mu_{\R,\q}^1)(\Psi_{M-B}^{-1}C_{\ba})}
=\frac{(\tilde\Phi_*\mu_{\R}^1)(\Psi_{M-B}^{-1}C_{\ba\hat
a\ba})}{(\tilde\Phi_*\mu_{\R}^1)(\Psi_{M-B}^{-1}C_{\ba})}.
\end{equation}

Observe that $\Psi_{M-B}^{-1}C_{\ba}=C\cap\Omega_\q^\mathbb Z$,
where $C$ is a cylinder in $\Omega_{\A,B}^\mathbb Z$ whose support
belongs to $\mathbb Z_+$. It follows that
\begin{equation*}
\psi^{-1}(\psi(C\cap\Omega_\q^\mathbb Z))=C\cap\Omega_\q^{\mathbb
Z_+},
\end{equation*}
where $\Omega_\q^{\mathbb Z_+}$ consists of all
$\omega=(\omega_k,k\in\mathbb Z)\in\Omega_{\A,B}^\mathbb Z$ such
that $(\omega_n,\dots,\omega_{n+|\q|-1})=\q$ for $n=0$ and for
infinitely many $n>0$. Hence
\begin{align*}
 (\psi_*(\tilde\Phi_*\mu_\R^1))(\psi
\Psi_{M-B}^{-1}C_\ba))&=(\psi_*(\tilde\Phi_*\mu_\R^1))(\psi
(C\cap\Omega_\q^\mathbb Z))\notag\\
&=({\tilde\Phi}_*\mu_\R^1)(C\cap\Omega_\q^{\mathbb Z_+}).
\end{align*}
Since the measure $\tilde\Phi_*\mu_{\R}^1$ is shift-invariant, while
the word $\q$ is a prefix of $\ba_1\in\wab$, from the Poincar{\'e}
Recurrence Theorem we obtain
$$
(\Phi_*\mu_\R^1)(C\cap\Omega_\q^\mathbb
Z)=(\Phi_*\mu_\R^1)(C\cap\Omega_\q^{\mathbb
Z_+})=(\Phi_*\mu_\R^1)(C),
$$
so that
\begin{equation}
\label{aaa} (\psi_*(\tilde\Phi_*\mu_\R^1))(\psi
\Psi_{M-B}^{-1}C_\ba))=({\tilde\Phi}_*\mu_\R^1)(\Psi_{M-B}^{-1}C_\ba).
\end{equation}

In this argument one may replace $C_\ba$ by $C_{\ba\hat a\ba}$ to
obtain
\begin{equation}
\label{aaa1}
(\psi_*(\tilde\Phi_*\mu_\R^1))(\psi\Psi_{M-B}^{-1}C_{\ba\hat
a\ba})=({\tilde\Phi}_*\mu_\R^1)(\Psi_{M-B}^{-1}C_{\ba\hat a\ba}).
\end{equation}
Substitution of (\ref{aaa}) and (\ref{aaa1}) in (\ref{ratio3}) with
taking into account (\ref{proj_nu}) yields
\begin{equation}
\label{ratio4} \frac{\bar\mu_\q(C_{\ba\hat a\ba})}
{\bar\mu_\q(C_\ba)}=\frac{(\Phi_*\nu)(\psi\Psi_{M-B}^{-1}C_{\ba\hat
a\ba})} {(\tilde\Phi_*\nu)(\psi\Psi_{M-B}^{-1}C_{\ba})}=
\frac{\nu(\Phi^{-1}(\psi\Psi_{M-B}^{-1}C_{\ba\hat
a\ba}))}{\nu(\Phi^{-1}(\psi\Psi_{M-B}^{-1}C_\ba))}.
\end{equation}

In a similar way as above we have
$$
\psi(\Psi_{M-B}^{-1}C_\ba)=C_{\ba(\q)}\cap\Omega_\q,\ \
\psi(\Psi_{M-B}^{-1}C_{\ba\hat a\ba})=C_{{\bf\hat
a(\q)}}\cap\Omega_\q,
$$
where $C_{\ba(\q)}$, $C_{{\bf\hat a(\q)}}$ are the cylinders in
$\Omega_{\A,B}$ corresponding to the words $\ba(\q):=\q q_1\q\dots\q
q_n\q\in\wab$ and ${\bf\hat a}(\q):=\ba(\q)\hat w\ba(\q)\in\wab$
with some $\hat w\in\Omega_{\A,B}$, respectively. Therefore (see
(\ref{ratio4}), (\ref{Delta_w3})),
\begin{equation}
\label{ratio5}
 \frac{\bar\mu_\q(C_{\ba\hat
a\ba})}{\bar\mu_\q(C_\ba)}=\frac{\nu(\Delta({\bf\hat
a}(\q)))}{\nu(\Delta(\ba(\q)))}.
\end{equation}

Let us now apply Proposition \ref{unifexp} with $\p'=\q q_1\q\dots\q
q_n$, $\p=\ba(\q)$, ${\bf r}=\hat{\bf a}(\q)$, where $n$ is large
enough. We may do so because of the following two facts: 1) the
choice of the word $\hat a$ above implies that $\hat w$ does not
contain subwords equal to $\ba(\q)$; 2) ${\rm d}(\p)\le
e^{-n\alpha}$ for some $\alpha>0$, as can be shown in the manner of
Subsection \ref{roof}. By this proposition combined with
(\ref{ratio3}) we obtain
$$
\left|\frac{\bar\mu_\q(C_{\ba\hat
a\ba})}{\bar\mu_\q(C_\ba)}-\exp(-m\tau_{\ba(\q)}^1(\la,\pi))
\right|\le\beta_2\exp(-n\alpha-m\tau_\q^1(\la,\pi)),
$$
where
$$
(\la,\pi)\in(\Phi^{-1}\psi\Psi_{M-B}C_{\ba\hat
a\ba})\cap\Phi^{-1}\Omega_{\ba(\q)}.
$$
Using (\ref{tau}), (\ref{roof1}), and (\ref{tau4})--(\ref{tau6}),
one can check that if $u\in C_{\ba\hat a\ba}$ is such that
$(u_n,\dots,u_{n+|\ba|-1})=\ba$ for infinitely many $n>0$ and $n<0$,
then $\tilde\tau(u,C_\ba)=\tau_{\ba(\q)}(\la,\pi)$, where
$(\la,\pi)=\Phi^{-1}\psi\Psi_{M-B}u\in\Phi^{-1}\Omega_{\ba(\q)}$. It
remains to note that by the Recurrence Theorem, $\bar\mu_\q$-almost
all $u\in C_\ba\hat a\ba$ satisfy the stated condition.

We thus see that the measure $\bar\mu_\q$ satisfies the assumptions
of Theorem \ref{main} (in particular, the constant $s$ in
(\ref{marg}) equals $m=h_{\Psi_*\mu_\R}(\sigma_\q,f_\q)$). This
theorem now implies that $(\bar\mu_\q)_f$ is a measure with maximal
entropy for $(\sigma_\q,f_\q)$ and hence $\mu_\R$ is a measure with
maximal entropy for $\{P^t\}$, so that $h_{\text
{top}}(\{P^t\})=h_{\text {top}}(\sigma_\q,f_\q)=m$.

As to the measure $\mu$, it follows that $h_\mu(\{P^t\})=
h_{\mu_\R}(\{P^t\})$. Then $\Psi_*\mu$ has the same entropy with
respect to the suspension flow $(\sigma_\q,f_\q)$ (see
(\ref{entropy1})) and hence is a measure with maximal entropy for
this flow. But from Theorem \ref{main} we know that such a measure
is unique. This completes the proof of Theorem \ref{maxentropy}.

{\bf Acknowledgements.}
A.I.B. is supported in part by Grant MK-4893.2010.1 of the President 
of the Russian Federation, by the Programme on Mathematical Control Theory of
the Presidium of the Russian Academy of Sciences,
by the Russian Foundation for Basic Research under
Grant 10-01-00739-à,
by the Russian Ministry of Education and Research
under the 2010 Programme for Development of Higher Education,
by the National
Science Foundation under grant DMS 0604386 and by the Edgar Odell
Lovett Fund at Rice University.
 B.M.G. is supported in part by the
Russian Foundation for Basic Research grants 07-01-92215 CNRS(L) and
08-01-00105.

\end{document}